\documentclass[british]{article}
\usepackage{ae,aecompl}
\usepackage[T1]{fontenc}
\usepackage[latin9]{inputenc}
\usepackage{enumitem}
\usepackage{amsmath}
\usepackage{amsthm}
\usepackage{amssymb}
\usepackage{geometry}
\usepackage{tablefootnote}

\makeatletter

\providecommand{\tabularnewline}{\\}

\theoremstyle{plain}
\newtheorem{thm}{\protect\theoremname}
\theoremstyle{definition}
\newtheorem{defn}[thm]{\protect\definitionname}
\theoremstyle{remark}
\newtheorem{rem}[thm]{\protect\remarkname}
\theoremstyle{plain}
\newtheorem{lem}[thm]{\protect\lemmaname}
\theoremstyle{plain}
\newtheorem{cor}[thm]{\protect\corollaryname}
\theoremstyle{plain}
\newtheorem{prop}[thm]{\protect\propositionname}
\theoremstyle{definition}
\newtheorem{example}[thm]{\protect\examplename}

\date{}
\usepackage{bbm}
\usepackage{array}
\newtheorem{notation}[thm]{Notation}
\makeatother

\usepackage{babel}
\providecommand{\corollaryname}{Corollary}
\providecommand{\definitionname}{Definition}
\providecommand{\examplename}{Example}
\providecommand{\lemmaname}{Lemma}
\providecommand{\propositionname}{Proposition}
\providecommand{\remarkname}{Remark}
\providecommand{\theoremname}{Theorem}

\begin{document}
\global\long\def\IN{\mathbb{N}}%
\global\long\def\II{\mathbbm{1}}%
\global\long\def\IZ{\mathbb{Z}}%
\global\long\def\IQ{\mathbb{Q}}%
\global\long\def\IR{\mathbb{R}}%
\global\long\def\IC{\mathbb{C}}%
\global\long\def\IP{\mathbb{P}}%
\global\long\def\IE{\mathbb{E}}%

\title{Probabilistic Voting Models\\
with Varying Speeds of Correlation Decay}
\author{Gabor Toth\thanks{IIMAS-UNAM, Mexico City, Mexico, gabor.toth@iimas.unam.mx }}
\maketitle
\begin{abstract}
\noindent We model voting behaviour in the multi-group setting of
a two-tier voting system using sequences of de Finetti measures. Our
model is defined by using the de Finetti representation of a probability
measure (i.e.\! as a mixture of conditionally independent probability
measures) describing voting behaviour. The de Finetti measure describes
the interaction between voters and possible outside influences on
them. We assume that for each population size there is a (potentially)
different de Finetti measure, and as the population grows, the sequence
of de Finetti measures converges weakly to the Dirac measure at the
origin, representing a tendency toward weakening social cohesion as
the population grows large. The resulting model covers a wide variety
of behaviour, ranging from independent voting in the limit under fast
convergence, a critical convergence speed with its own pattern of
behaviour, to a subcritical convergence speed which yields a model
in line with empirical evidence of real-world voting data, contrary
to previous probabilistic models used in the study of voting. These
models can be used, e.g., to study the problem of optimal voting weights
in two-tier voting systems.
\end{abstract}
Keywords: probabilistic voting models, spin models, de Finetti representation,
limit theorems, phase transitions

2020 Mathematics Subject Classification: 60F05, 82B20, 91B12

\section{Introduction}

Spin models from statistical mechanics have been used to model voting
behaviour (see e.g.\! \cite{BD,CGh2007,SRL,KL1,OEA,Toth,LSV}). These
stochastic models of voting describe the typical voting behaviour
in terms of a probability distribution over the set of all possible
voting outcomes for a population of $n\in\IN$ voters that decide
on a binary matter. A complete record of all individual decisions
is called a voting configuration $\left(x_{1},\ldots,x_{n}\right)\in\left\{ -1,1\right\} ^{n}$.
These binary choices, with values encoded as $\pm1$, cover both referenda
and elections for public office with two candidates. The individual
votes are random variables $X_{i}$, $i=1,2,\ldots,n$, with joint
distribution $\IP$. In the literature, it is customary to impose
the symmetry condition
\begin{equation}
\IP\left(X_{1}=x_{1},\ldots,X_{n}=x_{n}\right)=\IP\left(X_{1}=-x_{1},\ldots,X_{n}=-x_{n}\right)\label{eq:symmetry}
\end{equation}
for all voting configurations $\left(x_{1},\ldots,x_{n}\right)$.
The justification for this is that there is no fundamental distinction
between the two alternatives. We assume that each voter has fixed
preferences concerning all possible issues that can be put to vote.
However, the way the question is phrased in the referendum, or the
order in which the candidates appear on the ballot, is assumed to
be random, so it is reasonable to expect any given voter to choose
$+1$ with the same probability as $-1$.
\begin{defn}
\label{def:voting_measure}A probability measure on $\left\{ -1,1\right\} ^{n}$
that satisfies (\ref{eq:symmetry}) is called a \emph{voting measure}.
\end{defn}

With only two alternatives, we can apply the majority rule to determine
the winning alternative, and measure the outcome of the vote by adding
all votes to obtain the voting margin $S_{n}:=\sum_{i=1}^{n}X_{i}$.
If the voting margin is positive, the decision is in favour; otherwise,
it is against the proposal. The expected per capita absolute voting
margin is $\IE\left(\left|S_{n}\right|/n\right)$, and it captures
the essence of voting behaviour under the model described by $\IP$.

We would like to note that even though the motivation for the study
of the models in this article is voting theory, all of these models
can also be considered spin models in the context of statistical physics,
where there is a tendency for the spins to align.

The focus of this article is the asymptotic behaviour of sequences
of voting measures as the population size $n$ goes to infinity. Empirical
data of elections give us an idea of the properties these models should
have to reflect real-world voting. Penrose's square root law \cite{Penrose}
has been the gold standard in the field of voting theory. Penrose
investigated the question of how to equalise the influence of each
voter on the decision in a two-tier voting system. A two-tier voting
system consists of several groups or constituencies of possibly different
sizes and a council that makes decisions on behalf of the union of
these groups. Each group sends a representative to the council who
votes according to the preferences of the group. The voting weights
in the council have to take into account the potentially different
group sizes. According to Penrose's square root law, the optimal weights
in the council should be proportional to the square root of the size
of each group. Aside from equalising the influence of each voter on
the council vote, other criteria have been proposed. This is the problem
of optimal voting weights, a topic treated by a variety of authors
\cite{SZ1,SZ2,SZZ,SZ3,SRL,KL1,SZ4,Toth,KT_CBM,KT_CW_opt_weights}.
It is but one of the applications for voting models with real-world
consequences.

A different but related subject is the study of dynamic models of
binary voting which can also be interpreted as models of opinion dynamics.
One work which shares some overlap with the present article is \cite{Toth22},
in which we specify a dynamic model of binary voting, where the population
is subdivided in a large number of groups (such as work groups or
friendship groups) that discuss the issue to be voted on and influence
each other. After discussion takes place, there is a general tendency
for everybody to adopt the majority opinion prevalent in the group.
However, we allow for contrarian tendencies which manifest in a refusal
to adopt said opinion. This is triggered randomly for each member
of the group and the probability is what we call a `flip parameter'.
Several rounds of these discussions take place sequentially, updating
opinions across the population in each round. This type of model is
very general in that a wide variety of behaviours can be observed,
depending on the choice of the flip parameters. Therefore, in \cite{Toth22},
we study how the general class of these models behaves if we suppose
the flip parameters are randomly chosen with a certain correlation
structure inherent in the joint distribution of all flip parameters.
Thus, the main differences to the present work lie in the dynamic
aspect of the model considered in \cite{Toth22}, the small size of
groups, where instead of holding the number of groups fixed, we instead
study groups of fixed sizes, as well as the modelling approach of
selecting the parameters at random. In the present article, following
the approach in the works cited in the last paragraph, we study static
models of voting with fixed parameters.

A basic assumption in \cite{Penrose} is the independence of all voters.
Intuitively speaking, the assumption that all voters are independent
seems unrealistic. Gelman et al.\! \cite{GKT,GKB} criticised the
square root law based on statistical evidence. $\IE\left(\left|S_{n}\right|/n\right)$
is a characteristic of the distribution of voting configurations under
a voting measure which can be statistically estimated. Gelman et al.\!
pointed out that under the assumption of independence of all voters,
which leads to a binomial distribution of the sum of all votes, $\IE\left(\left|S_{n}\right|/n\right)$
should be of order $1/\sqrt{n}$ as $n$ goes to infinity. This would
imply that larger countries or constituencies should on average have
far closer elections, i.e.\!  smaller per capita voting margins,
than smaller countries. Gelman et al.\! argued that although this
can be observed in the data on U.S. and European elections to some
extent, the per capita voting margins decrease as a far lower power
of $n$ than the power $1/2$ specified by the square root law. Indeed,
the typical magnitude of $\IE\left(\left|S_{n}\right|/n\right)$ seems
to be of order $n^{-\alpha}$ with $0.1\leq\alpha\leq0.2$. This has
implications for the determination of the optimal voting weights in
two-tier voting systems.

Models of voting behaviour are often adapted from statistical physics.
Two categories are the collective bias model and the Curie-Weiss model.
See Section \ref{sec:Models} for a discussion of these models and
how they fit into the framework introduced in Definition \ref{def:de_Finetti_model}
below. 

In this article, we will deal with the multivariate (i.e.\!  multi-group)
setting. This is due to the application of these models to two-tier
voting systems, in which the overall population is subdivided into
a fixed number of groups such as e.g.\! states of a federal republic.
We could describe each group's voters by a separate voting model,
but this would preclude us from considering correlated voting across
group boundaries. Therefore, we choose to analyse multi-group models.
The general framework is a set of binary random variables double-indexed
by the group and individual $X_{\lambda i}\in\left\{ -1,1\right\} $,
where $\lambda=1,\ldots,M$ indicates the group and $i=1,\ldots,n_{\lambda}$,
the individual voter. The group sizes $n_{\lambda}$ sum to $n$,
the overall population size. We will assume that all $M$ groups grow
without bound as the overall population goes to infinity. Each $n_{\lambda}$
is thus actually a sequence $\left(n_{\lambda}(n)\right)_{n\in\IN}$
that goes to infinity as $n\rightarrow\infty$. So for each value
of the overall population $n$, we have a vector of group sizes $\left(n_{1},\ldots,n_{M}\right)$.

Our aim in this article is to study sequences of de Finetti measures
$\left(\mu_{n}\right)_{n}$ with certain concentration properties.
Aside from the case of each $\mu_{n}$ living on the compact set $\left[-1,1\right]^{M}$,
we will also consider more general measures with supports in $\IR^{M}$.
All these sequences will have in common that they converge weakly
to $\delta_{0}$, the Dirac measure at the origin.  This convergence
assumption represents a tendency towards lower social cohesion or
lower central influence as the size of the population increases.

We will use the Rademacher distribution in the definition of de Finetti
voting models below.
\begin{defn}
\label{def:Rademacher}The \emph{Rademacher distribution} with parameter
$s\in\left[-1,1\right]$ is a probability measure $P_{s}$ on $\left\{ -1,1\right\} $
given by $P_{s}\left\{ 1\right\} :=\frac{1+s}{2}$.
\end{defn}

\begin{notation}\label{not:prod_measure}We will write for all $t\in\left[-1,1\right]^{M}$
and all $\left(x_{11},\ldots,x_{Mn_{M}}\right)\in\left\{ -1,1\right\} ^{n}$
\[
P_{t}^{\otimes n}\left(x_{11},\ldots,x_{Mn_{M}}\right):=\prod_{\lambda=1}^{M}\prod_{i=1}^{n_{\lambda}}P_{t_{\lambda}}\left\{ x_{\lambda i}\right\} .
\]
Let $E_{t}$ stand for the expectation with respect to the product
measure $P_{t}^{\otimes n}$.\end{notation}

Recall that $M\in\IN$ is the number of groups and $n=\sum_{\lambda=1}^{M}n_{\lambda}\in\IN$
is the size of the overall population. 
\begin{defn}
\label{def:de_Finetti_model}Let $\left(\mu_{n}\right)_{n}$ be a
sequence of probability measures on $\IR^{M}$ which are symmetric
with respect to the origin, i.e. $\mu_{n}\left(A\right)=\mu_{n}\left(-A\right)$
holds for all measurable sets $A\subset\IR^{M}$ and all $n\in\IN$.
The voting measure $\IP_{n}$ is defined by
\[
\IP_{n}\left(X_{11}=x_{11},\ldots,X_{Mn_{M}}=x_{Mn_{M}}\right):=\int_{\IR^{M}}P_{\bar{m}}^{\otimes n}\left(x_{11},\ldots,x_{Mn_{M}}\right)\mu_{n}\left(\textup{d}m\right)
\]
for all voting configurations $\left(x_{11},\ldots,x_{Mn_{M}}\right)\in\left\{ -1,1\right\} ^{n}$
and all $n\in\IN$. The expression $\bar{m}$ is some function $m\in\IR^{M}\mapsto\bar{m}\in\IR^{M}$
with the following properties:
\begin{itemize}
\item $m\mapsto\bar{m}$ is increasing in each component of $m$ and its
range lies in $\left[-1,1\right]^{M}$.
\item For each $\lambda$, $\lim_{m\rightarrow0}\frac{\bar{m}_{\lambda}}{m_{\lambda}}=1$.
\item The limit of $\bar{m}$ taken where each component of $m$ goes to
$\infty$ is 1, and the limit where each component goes to $-\infty$
is $-1$.
\end{itemize}
We will call such sequences of voting measures $\left(\IP_{n}\right)_{n}$
\emph{de Finetti voting models} with de Finetti sequence $\left(\mu_{n}\right)_{n}$.
We will refer to $\left(\mu_{n}\right)_{n}$ as the sequence of de
Finetti measures of $\left(\IP_{n}\right)_{n}$.
\end{defn}

\begin{rem}
Note that the assumed symmetry of each probability measure $\mu_{n}$
implies the symmetry condition (\ref{eq:symmetry}) for each $\IP_{n}$.
\end{rem}

\begin{rem}
We will concentrate on the case where $\left(\mu_{n}\right)_{n}$
is supported on $\left[-1,1\right]^{M}$ with $\bar{m}=m$. The reason
we consider the more general setup with supports not contained in
$\left[-1,1\right]^{M}$ is that some voting models such as the Curie-Weiss
model defined in Section \ref{sec:Models} fall into this category.
However, there are no significant mathematical differences between
the more restrictive assumption of a compact support and the general
case.
\end{rem}

For large population sizes $n$, Definition \ref{def:de_Finetti_model}
covers all models with exchangeable random variables $X_{\lambda i}$
within each group $\lambda$. (See \cite{DF1980} for the result that
states that for finite numbers $n$ of exchangeable random variables,
the total variation distance between the joint distribution of the
random variables and the distribution of a mixture as in the formula
above is at most of order $1/n$.) In the language of voting theory,
these voting systems are called neutral, as the identity of the individual
voters (at least within each group) does not affect the outcome.
\begin{defn}
\label{def:voting_margins}Suppose $\left(\IP_{n}\right)_{n}$ is
a de Finetti voting model. We define the \emph{group voting margins}
for each group $\lambda$:
\[
S_{n,\lambda}:=\sum_{i=1}^{n_{\lambda}}X_{\lambda i}.
\]
\end{defn}

Instead of the voting margin $S_{n}$ of a single-group model, we
will be interested in the joint distribution of the group voting margins
$S_{n,1},\ldots,S_{n,M}$; more precisely, we will study the asymptotic
behaviour of the sums
\begin{equation}
\boldsymbol{S}_{n}:=\left(\frac{S_{n,1}}{\gamma_{n,1}},\ldots,\frac{S_{n,M}}{\gamma_{n,M}}\right),\label{eq:S^n}
\end{equation}
normalised by sequences $\left(\gamma_{n,\lambda}\right)_{n\in\IN}$
to ensure $\boldsymbol{S}_{n}$ converges to some limiting distribution
as $n$ goes to infinity. 

We will stipulate that $\left(\mu_{n}\right)_{n}$ converges weakly
to $\delta_{0}$. This represents an underlying tendency towards lower
social cohesion which is reflected in the convergence to 0 of the
correlation $\IE\left(X_{\lambda1}X_{\nu2}\right)$ between two votes
belonging any groups $\lambda$ and $\nu$. Note that this does not
necessarily imply that the limiting distribution of $\boldsymbol{S}_{n}$
features independent entries. We will return to the topic of the correlations
$\IE\left(X_{\lambda1}X_{\nu2}\right)$ in Proposition \ref{prop:corr}
and in Section \ref{sec:Models}.

In Section \ref{sec:Results}, we present limit theorems for $\boldsymbol{S}_{n}$
given sequences $\left(\mu_{n}\right)_{n}$ that converge weakly to
$\delta_{0}$ and show that the speed of convergence is crucial for
the limiting distribution. The main results are limit Theorems \ref{thm:fast_conv}
and \ref{thm:phase_transition} and Corollary \ref{cor:voting_margins}
concerning the behaviour of the expected per capita voting margins.
Section \ref{sec:Proofs} contains the proofs of these results. Section
\ref{sec:Models} discusses some voting models featured prominently
in the literature, how they fit into the framework introduced in Definition
\ref{def:de_Finetti_model}, and in Example \ref{exa:de_Finetti}
provides a de Finetti voting model which fits the empirical evidence.
Finally, Section \ref{sec:Conclusion} presents the conclusions.

\section{\label{sec:Results}Results}

\begin{notation}In the following, let for any vectors $x,y\in\IR^{d}$
the expression $x\circ y\in\IR^{d}$ stand for the componentwise multiplication
of $x$ and $y$. We will use the symbol `$\ast$' to denote the
convolution of two measures. `$\xrightarrow[n\rightarrow\infty]{\textup{d}}$'
stands for convergence in distribution as $n$ goes to infinity. We
will write $\mathcal{N}\left(0,C\right)$ for a centred multivariate
normal distribution with covariance matrix $C$, and $\delta_{x}$
for the Dirac measure at $x\in\IR^{d}$. $I_{k}$ will stand for the
identity matrix of dimension $k\in\IN$. For $a,b\in\IR^{M}$ with
$a\leq b$ componentwise, we will use $\left[a,b\right]$ to denote
$\prod_{\lambda=1}^{M}\left[a_{\lambda},b_{\lambda}\right]$. Let
$\left\Vert \cdot\right\Vert _{\infty}$ stand for the sup norm on
$\IR^{d}$ for any $d\in\IN$ as well as on the set of all functions
$f:\IR^{d}\rightarrow\IR$.\end{notation}

We will be working with sequences of vectors $\left(\varepsilon_{n}\right)_{n}$
in $\left(0,\infty\right)^{M}$ such that $\left\Vert \varepsilon_{n}\right\Vert _{\infty}\xrightarrow[n\rightarrow\infty]{}0$.
\begin{lem}
\label{lem:weak_conv}Let $\left(\mu_{n}\right)_{n}$ be a sequence
of probability measures on $\IR^{M}$. Then the following statements
are equivalent:
\begin{enumerate}[label=(\roman*)]
\item $\left(\mu_{n}\right)_{n}$ converges weakly to $\delta_{0}$.
\item For all $\delta>0$, we have $\mu_{n}\left(\prod_{\lambda=1}^{M}\left[-\delta,\delta\right]\right)\xrightarrow[n\rightarrow\infty]{}1$.
\item There is a sequence $\left(\varepsilon_{n}\right)_{n}$\textup{ in
$\left(0,\infty\right)^{M}$ with} $\left\Vert \varepsilon_{n}\right\Vert _{\infty}\xrightarrow[n\rightarrow\infty]{}0$
\textup{such that $\mu_{n}\left(\left[-\varepsilon_{n},\varepsilon_{n}\right]\right)\xrightarrow[n\rightarrow\infty]{}1$}.
\end{enumerate}
\end{lem}

This lemma states that by considering all sequences $\left(\varepsilon_{n}\right)_{n}$
as in $(iii)$ above, we are covering all sequences of probability
measures $\left(\mu_{n}\right)_{n}$ that converge weakly to $\delta_{0}$.

We will use the following notation for the asymptotic behaviour of
sequences:

\begin{notation}For any complex-valued sequences $\left(f_{n}\right)_{n},\left(g_{n}\right)_{n}$,
we will write $f_{n}=\Theta\left(g_{n}\right)$ if
\[
0<\liminf_{n\rightarrow\infty}\left|\frac{f_{n}}{g_{n}}\right|\leq\limsup_{n\rightarrow\infty}\left|\frac{f_{n}}{g_{n}}\right|<\infty,
\]
$f_{n}=O\left(g_{n}\right)$ if
\[
\limsup_{n\rightarrow\infty}\left|\frac{f_{n}}{g_{n}}\right|<\infty,
\]
and $f_{n}=o\left(g_{n}\right)$ if
\[
\lim_{n\rightarrow\infty}\frac{f_{n}}{g_{n}}=0.
\]
\end{notation}

The convergence speed of $\left(\mu_{n}\right)_{n}$ turns out to
be crucial for the limiting distribution of the suitably normalised
sums $\boldsymbol{S}_{n}$. First we present a theorem that says that
under `fast convergence', a term we will give formal meaning to
in Definition \ref{def:conv_speed}, we obtain a universal limit for
the normalised sums which is independent of the specific sequence
of de Finetti measures, with only the convergence speed mattering.
\begin{thm}
\label{thm:fast_conv}Let $\left(\IP_{n}\right)_{n}$ be a de Finetti
voting model with group sizes $n_{\lambda}$, group voting margins
$S_{n,\lambda}$, $\lambda=1,\ldots,M$, and de Finetti measures $\left(\mu_{n}\right)_{n}$.
Let $\left(\varepsilon_{n}\right)_{n}$ be a sequence in $\left(0,\infty\right)^{M}$
with $\varepsilon_{n,\lambda}=o\left(\frac{1}{\sqrt{n_{\lambda}}}\right)$
for each coordinate $\lambda$. If $\mu_{n}\left(\left[-\varepsilon_{n},\varepsilon_{n}\right]\right)\xrightarrow[n\rightarrow\infty]{}1$,
then
\[
\left(\frac{S_{n,1}}{\sqrt{n_{1}}},\ldots,\frac{S_{n,M}}{\sqrt{n_{M}}}\right)\xrightarrow[n\rightarrow\infty]{\textup{d}}\mathcal{N}\left(0,I_{M}\right).
\]
\end{thm}

This theorem says that if convergence is fast enough -- that is,
faster than $1/\sqrt{n_{\lambda}}$ in each component -- then the
specifics of the sequence of distributions $\left(\mu_{n}\right)_{n}$
are lost in the limit. The groups become independent in the limit,
and each normalised group voting margin $S_{n,\lambda}/\sqrt{n_{\lambda}}$
converges to a standard normal distribution. As we will next see,
this universality is lost if the convergence speed is not faster than
the critical speed $1/\sqrt{n_{\lambda}}$. Then, the specifics of
the sequence $\left(\mu_{n}\right)_{n}$ affect the limiting distribution
of the normalised sums. For this reason, we assume that each $\mu_{n}$
is a contraction of some fixed probability measure $\mu$ on $\IR^{M}$.

For the next theorem, we define the notation for a contraction:

\begin{notation}If $\mu$ is a probability measure on $\IR^{d}$
and $h\in\left(0,\infty\right)^{d}$, then we denote by $h\circ\mu$
the rescaled measure of $\mu$ by $h$ given by $h\circ\mu\left(A\right):=\mu\left(h\circ A\right)$
for all measurable sets $A\subset\IR^{d}$.\end{notation}
\begin{thm}
\label{thm:phase_transition}Let $\mu$ be a probability measure on
$\IR^{M}$ and $\left(\varepsilon_{n}\right)_{n}$ a sequence in $\left(0,\infty\right)^{M}$
with $\left\Vert \varepsilon_{n}\right\Vert _{\infty}\xrightarrow[n\rightarrow\infty]{}0$.
We define for each $n\in\IN$ the probability measure $\mu_{n}:=\left(\frac{1}{\varepsilon_{n,1}},\ldots,\frac{1}{\varepsilon_{n,M}}\right)\circ\mu$.
Let $\left(\IP_{n}\right)_{n}$ be a de Finetti voting model with
group sizes $n_{\lambda}$, group voting margins $S_{n,\lambda}$,
$\lambda=1,\ldots,M$, and de Finetti measures $\left(\mu_{n}\right)_{n}$.
\begin{enumerate}
\item If for each component $\lambda$ $\varepsilon_{n,\lambda}=o\left(\frac{1}{\sqrt{n_{\lambda}}}\right)$,
then
\[
\left(\frac{S_{n,1}}{\sqrt{n_{1}}},\ldots,\frac{S_{n,M}}{\sqrt{n_{M}}}\right)\xrightarrow[n\rightarrow\infty]{\textup{d}}\mathcal{N}\left(0,I_{M}\right).
\]
\item If for each component $\lambda$ $\lim_{n\rightarrow\infty}\varepsilon_{n,\lambda}\sqrt{n_{\lambda}}=h_{\lambda}>0$,
we set $h:=\left(h_{1},\ldots,h_{M}\right)$, and then
\[
\left(\frac{S_{n,1}}{\sqrt{n_{1}}},\ldots,\frac{S_{n,M}}{\sqrt{n_{M}}}\right)\xrightarrow[n\rightarrow\infty]{\textup{d}}\mathcal{N}\left(0,I_{M}\right)\ast\left(h\circ\mu\right).
\]
\item If for each component $\lambda$ $\frac{1}{\sqrt{n_{\lambda}}}=o\left(\varepsilon_{n,\lambda}\right)$,
we set for each $n$ and each $\lambda$ $\gamma_{n,\lambda}:=\varepsilon_{n,\lambda}n_{\lambda}$,
and then
\[
\left(\frac{S_{n,1}}{\gamma_{n,1}},\ldots,\frac{S_{n,M}}{\gamma_{n,M}}\right)\xrightarrow[n\rightarrow\infty]{\textup{d}}\mu.
\]
\end{enumerate}
\end{thm}

Under fast contraction, i.e.\!  $\varepsilon_{n,\lambda}=o\left(1/\sqrt{n_{\lambda}}\right)$
for each $\lambda$, we obtain the universal limiting distribution
$\mathcal{N}\left(0,I_{M}\right)$ as in Theorem \ref{thm:fast_conv}.
However, when the contraction speed is critical, i.e.\!  $\lim_{n\rightarrow\infty}\varepsilon_{n,\lambda}\sqrt{n_{\lambda}}=h_{\lambda}>0$
for each $\lambda$, then the limiting distribution for the normalised
sums is the convolution of the rescaled version $h\circ\mu$ of the
measure $\mu$, which is used to define the contracting sequence $\left(\mu_{n}\right)_{n}$,
and a multivariate normal noise with independent components and standard
normal marginal distributions. Lastly, if the contraction speed is
subcritical, then the sums must be normalised differently than in
the supercritical and critical cases. The normalisation factor is
$\gamma_{n,\lambda}=\varepsilon_{n,\lambda}n_{\lambda}$ for each
component $\lambda$. The assumption $1/\sqrt{n_{\lambda}}=o\left(\varepsilon_{n,\lambda}\right)$
implies that $\gamma_{n,\lambda}$ goes to infinity faster than $\sqrt{n_{\lambda}}$.
The limiting distribution of the normalised sums in the subcritical
case is the measure $\mu$.

From Theorem \ref{thm:phase_transition}, the corresponding result
for single-group models (with $M=1$ and $n_{1}=n$) follows:
\begin{cor}
\label{cor:single_group}Let $\mu$ be a probability measure on $\IR$
and $\left(\varepsilon_{n}\right)_{n}$ a sequence in $\left(0,\infty\right)$
with $\varepsilon_{n}\xrightarrow[n\rightarrow\infty]{}0$. We define
for each $n\in\IN$ the probability measure $\mu_{n}:=\frac{1}{\varepsilon_{n}}\circ\mu$.
Let $\left(\IP_{n}\right)_{n}$ be a de Finetti voting model with
de Finetti measures $\left(\mu_{n}\right)_{n}$.
\begin{enumerate}
\item If $\varepsilon_{n}=o\left(\frac{1}{\sqrt{n}}\right)$, then
\[
\frac{S_{n}}{\sqrt{n}}\xrightarrow[n\rightarrow\infty]{\textup{d}}\mathcal{N}\left(0,1\right).
\]
\item If $\lim_{n\rightarrow\infty}\varepsilon_{n}\sqrt{n}=h>0$, then
\[
\frac{S_{n}}{\sqrt{n}}\xrightarrow[n\rightarrow\infty]{\textup{d}}\mathcal{N}\left(0,1\right)\ast\left(h\circ\mu\right).
\]
\item If $\frac{1}{\sqrt{n_{\lambda}}}=o\left(\varepsilon_{n}\right)$,
we set for each $n$ $\gamma_{n}:=\varepsilon_{n}n$. Then
\[
\frac{S_{n}}{\gamma_{n}}\xrightarrow[n\rightarrow\infty]{\textup{d}}\mu.
\]
\end{enumerate}
\end{cor}

The next definition formalises the different convergence speeds in
Theorem \ref{thm:phase_transition}. We will classify some voting
models that fit the de Finetti framework in Section \ref{sec:Models}
according to the speed of convergence to $\delta_{0}$ of their de
Finetti measures $\left(\mu_{n}\right)_{n}$.
\begin{defn}
\label{def:conv_speed}Let $\left(\IP_{n}\right)_{n}$ be a de Finetti
voting model with group sizes $n_{\lambda}$, $\lambda=1,\ldots,M$,
de Finetti measures $\left(\mu_{n}\right)_{n}$, and $\left(\varepsilon_{n}\right)_{n}$
a sequence in $\left(0,\infty\right)^{M}$ with $\left\Vert \varepsilon_{n}\right\Vert _{\infty}\xrightarrow[n\rightarrow\infty]{}0$.
Assume $\mu_{n}\left(\left[-\varepsilon_{n},\varepsilon_{n}\right]\right)\xrightarrow[n\rightarrow\infty]{}1$
holds.
\begin{enumerate}
\item If for each component $\lambda$ $\varepsilon_{n,\lambda}=o\left(\frac{1}{\sqrt{n_{\lambda}}}\right)$,
we will say the model exhibits \emph{fast }(or \emph{supercritical})\emph{
convergence}.
\item If for each component $\lambda$ $\lim_{n\rightarrow\infty}\varepsilon_{n,\lambda}\sqrt{n_{\lambda}}=h_{\lambda}>0$,
we will say the model has \emph{critical convergence speed}.
\item If for each component $\lambda$ $\frac{1}{\sqrt{n_{\lambda}}}=o\left(\varepsilon_{n,\lambda}\right)$,
we will say the model exhibits \emph{slow }(or \emph{subcritical})
\emph{convergence}.
\end{enumerate}
\end{defn}

Theorem \ref{thm:phase_transition} states that there is a phase transition
in the space of convergence speeds with a critical speed of order
$1/\sqrt{n_{\lambda}}$, as well as supercritical and subcritical
regimes.

Next we state that different components (or groups) $\lambda$ can
be in different regimes depending on the contraction speed $\varepsilon_{n,\lambda}$
for each group.
\begin{thm}
\label{thm:clusters}Let $\left(\IP_{n}\right)_{n}$ be a de Finetti
voting model with group sizes $n_{\lambda}$, group voting margins
$S_{n,\lambda}$, $\lambda=1,\ldots,M$, de Finetti measures $\left(\mu_{n}\right)_{n}$
defined as a contraction of a probability measure $\mu$ on $\IR^{M}$
as in Theorem \ref{thm:phase_transition}, and $\left(\varepsilon_{n}\right)_{n}$
a sequence in $\left(0,\infty\right)^{M}$ with $\left\Vert \varepsilon_{n}\right\Vert _{\infty}\xrightarrow[n\rightarrow\infty]{}0$.
Let the $M$ groups be partitioned into three clusters $C_{1},C_{2},$
and $C_{3}$, each of which comprises $M_{i}>0$ of the $M$ groups,
respectively. Let the contraction speeds of $\left(\varepsilon_{n}\right)_{n}$
and the definition of the sequence $\left(\gamma_{n}\right)_{n}$
be those given in the following display:
\begin{align*}
\varepsilon_{n,\lambda} & =o\left(\frac{1}{\sqrt{n_{\lambda}}}\right), & \gamma_{n,\lambda} & :=\sqrt{n_{\lambda}}, & \lambda & \in C_{1},\\
\lim_{n\rightarrow\infty}\varepsilon_{n,\lambda}\sqrt{n_{\lambda}} & =h_{\lambda}>0, & \gamma_{n,\lambda} & :=\sqrt{n_{\lambda}}, & \lambda & \in C_{2},\\
\frac{1}{\sqrt{n_{\lambda}}} & =o\left(\varepsilon_{n,\lambda}\right), & \gamma_{n,\lambda} & :=\varepsilon_{n,\lambda}n_{\lambda}, & \lambda & \in C_{3}
\end{align*}

We set $h:=\left(h_{M_{1}+1},\ldots,h_{M_{1}+M_{2}},1,\ldots,1\right)\in\IR^{M_{2}+M_{3}}$.

Then
\[
\left(\frac{S_{n,1}}{\gamma_{n,1}},\ldots,\frac{S_{n,M}}{\gamma_{n,M}}\right)\xrightarrow[n\rightarrow\infty]{\textup{d}}\left(\mathcal{N}\left(0,I_{M_{1}}\right),\left(\mathcal{N}\left(0,I_{M_{2}}\right),\delta_{0}\right)\ast\left(h\circ\mu_{C_{2},C_{3}}\right)\right),
\]
where $\mu_{C_{2},C_{3}}$ is the marginal distribution of the coordinates
belonging to $C_{2}$ and $C_{3}$ under the measure $\mu$.
\end{thm}

Note that the groups in $C_{1}$ are asymptotically independent both
of each other as well as of the groups belonging to $C_{2}$ and $C_{3}$,
whereas the groups in $C_{2}$ and $C_{3}$ preserve the dependence
given by the probability measure $\mu$ with some additive noise introduced
in $C_{2}$ due to the critical nature of the contraction speed for
that cluster.

We turn to the behaviour of the expected per capita voting margins
$\IE\left(\left|S_{n,\lambda}\right|/n_{\lambda}\right)$, as these
can be used to assess the goodness of fit of the model to the data
presented in \cite{GKT,GKB}.
\begin{cor}
\label{cor:voting_margins}Let $\left(\IP_{n}\right)_{n}$ be a de
Finetti voting model with group sizes $n_{\lambda}$, group voting
margins $S_{n,\lambda}$, $\lambda=1,\ldots,M$, de Finetti measures
$\left(\mu_{n}\right)_{n}$ defined as a contraction of a probability
measure $\mu$ on $\IR^{M}$ as in Theorem \ref{thm:phase_transition},
and $\left(\varepsilon_{n}\right)_{n}$ a sequence in $\left(0,\infty\right)^{M}$
with $\left\Vert \varepsilon_{n}\right\Vert _{\infty}\xrightarrow[n\rightarrow\infty]{}0$.
\begin{enumerate}
\item If for some $\lambda$ $\varepsilon_{n,\lambda}=o\left(\frac{1}{\sqrt{n_{\lambda}}}\right)$,
then
\[
\IE\left(\frac{\left|S_{n,\lambda}\right|}{n_{\lambda}}\right)=\Theta\left(n_{\lambda}^{-1/2}\right).
\]
\item If for some $\lambda$ $\lim_{n\rightarrow\infty}\varepsilon_{n,\lambda}\sqrt{n_{\lambda}}=h_{\lambda}>0$,
then
\[
\IE\left(\frac{\left|S_{n,\lambda}\right|}{n_{\lambda}}\right)=\Theta\left(n_{\lambda}^{-1/2}\right).
\]
\item If for some $\lambda$ $\frac{1}{\sqrt{n_{\lambda}}}=o\left(\varepsilon_{n,\lambda}\right)$,
then
\[
\IE\left(\frac{\left|S_{n,\lambda}\right|}{n_{\lambda}}\right)=\Theta\left(\varepsilon_{n,\lambda}\right).
\]
\end{enumerate}
\end{cor}

Hence, if $\varepsilon_{n,\lambda}=n_{\lambda}^{-\alpha}$ with $0.1\leq\alpha\leq0.2$,
then we land in the empirical range given in \cite{GKT,GKB} for the
group voting margins $\IE\left(\left|S_{n,\lambda}\right|/n_{\lambda}\right)$.
In accordance with Definition \ref{def:conv_speed}, the empirical
evidence suggests that de Finetti voting models with \emph{subcritical}
convergence speeds best describe real-world voting behaviour. So it
seems that the assumption of a decay of social cohesion as the population
grows is correct; however, this decay is fairly slow. See Example
\ref{exa:de_Finetti} in Section \ref{sec:Models} for a simple example
of such a de Finetti voting model where the parameters can be adjusted
to fit the empirical per capita voting margins presented in \cite{GKT,GKB}.

Another way to assess the asymptotic loss of social cohesion is to
look at the behaviour of correlations $\IE\left(X_{\lambda1}X_{\nu2}\right)$
between two votes for any groups $\lambda$ and $\nu$. Note that
due to the exchangeability of the random variables $\left\{ X_{\lambda i}\,|\,i=1,\ldots,n_{\lambda}\right\} $
for any group $\lambda$, it does not matter which two votes we pick
for our correlation as long as they are two different votes:
\begin{align*}
\IE\left(X_{\lambda i}X_{\lambda j}\right) & =\IE\left(X_{\lambda i'}X_{\lambda j'}\right),\quad i,i',j,j'=1,\ldots,n_{\lambda},i\neq j,i'\neq j',\\
\IE\left(X_{\lambda i}X_{\nu j}\right) & =\IE\left(X_{\lambda i'}X_{\nu j'}\right),\quad\lambda\neq\nu,i,i'=1,\ldots,n_{\lambda},j,j'=1,\ldots,n_{\nu}.
\end{align*}
Recall the transformation $m\mapsto\bar{m}$ from Definition \ref{def:de_Finetti_model}.
\begin{prop}
\label{prop:corr}Let $\left(\IP_{n}\right)_{n}$ be a de Finetti
voting model and $\left(\varepsilon_{n}\right)_{n}$ a sequence in
$\left(0,\infty\right)^{M}$ with $\left\Vert \varepsilon_{n}\right\Vert _{\infty}\xrightarrow[n\rightarrow\infty]{}0$.
If $\mu_{n}\left(\left[-\varepsilon_{n},\varepsilon_{n}\right]\right)\xrightarrow[n\rightarrow\infty]{}1$,
then
\[
\IE\left(X_{\lambda1}X_{\nu2}\right)=\int_{\IR^{M}}\bar{m}_{\lambda}\bar{m}_{\nu}\,\mu_{n}\left(\textup{d}m\right)\xrightarrow[n\rightarrow\infty]{}0,\quad\lambda,\nu=1,\ldots,M.
\]
\end{prop}

The proposition says that as a consequence of the loss of social cohesion
represented by the assumption $\left\Vert \varepsilon_{n}\right\Vert _{\infty}\xrightarrow[n\rightarrow\infty]{}0$
and $\mu_{n}\left(\left[-\varepsilon_{n},\varepsilon_{n}\right]\right)\xrightarrow[n\rightarrow\infty]{}1$,
individual votes become asymptotically uncorrelated, even if both
votes belong to the same group. This is an interesting contrast to
the results above which state that, provided that convergence is \emph{not
supercritical}, the group voting margins do not, in general, become
independent in the large population limit.
\begin{rem}
Note that we cannot say anything about the speed of convergence of
$\IE\left(X_{\lambda1}X_{\nu2}\right)\xrightarrow[n\rightarrow\infty]{}0$
unless in addition to the speed of convergence of each $\varepsilon_{n,\lambda}\xrightarrow[n\rightarrow\infty]{}0$
we also make an assumption about the speed of convergence of $\mu_{n}\left(\left[-\varepsilon_{n},\varepsilon_{n}\right]\right)\xrightarrow[n\rightarrow\infty]{}1$.
\end{rem}

Finally, we give a local limit theorem for the fast convergence case.
In this regime, the convergence in distribution can be strengthened
to local convergence, i.e.\!  convergence of the scaled point probabilities
for the normalised sums to the Lebesgue density function $\phi$ of
$\mathcal{N}\left(0,I_{M}\right)$. For this, we need to make an assumption
about the tails of $\left(\mu_{n}\right)_{n}$.
\begin{thm}
\label{thm:LLT}Let $\left(\IP_{n}\right)_{n}$ be a de Finetti voting
model with group sizes $n_{\lambda}$, group voting margins $S_{n,\lambda}$,
$\lambda=1,\ldots,M$, and de Finetti measures $\left(\mu_{n}\right)_{n}$.
Let $\left(\varepsilon_{n}\right)_{n}$ be a sequence in $\left(0,\infty\right)^{M}$
with the property that for each component $\lambda$ $\varepsilon_{n,\lambda}=o\left(\frac{1}{\sqrt{n_{\lambda}}}\right)$.
If $\mu_{n}\left(\left[-\varepsilon_{n},\varepsilon_{n}\right]\right)\xrightarrow[n\rightarrow\infty]{}1$
and there is a constant $\tau\in\left(0,1\right)^{M}$ such that
\begin{equation}
\sum_{n=1}^{\infty}n^{\frac{M-1}{2}}\mu_{n}\left(\IR^{M}\backslash\left[-\tau,\tau\right]\right)<\infty,\label{eq:summability}
\end{equation}
then
\[
\sup_{x\in\mathcal{L}_{n}}\left|\frac{\prod_{\lambda=1}^{M}\sqrt{n_{\lambda}}}{2^{M}}\mathbb{P}\left(\left(\frac{S_{n,1}}{\sqrt{n_{1}}},\ldots,\frac{S_{n,M}}{\sqrt{n_{M}}}\right)=x\right)-\phi(x)\right|\xrightarrow[n\to\infty]{}0,
\]
where $\mathcal{L}_{n}$ is the lattice on which the vector $\left(\frac{S_{n,1}}{\sqrt{n_{1}}},\ldots,\frac{S_{n,M}}{\sqrt{n_{M}}}\right)$
lives.
\end{thm}

Note that the summability condition (\ref{eq:summability}) of the
sequence $\left(n^{\frac{M-1}{2}}\mu_{n}\left(\IR^{M}\backslash\left[-\tau,\tau\right]\right)\right)_{n}$
is a weaker assumption regarding the convergence speed of $\mu_{n}\left(\left[-\varepsilon_{n},\varepsilon_{n}\right]\right)\xrightarrow[n\rightarrow\infty]{}1$
than an exponential concentration property\footnote{The Curie-Weiss model defined in Section \ref{sec:Models} satisfies
such an exponential concentration condition (\ref{eq:exp_concentr_CW})
in its high temperature regime.} in which for each $\delta>0$ there are $C,D>0$ such that
\[
\mu_{n}\left(\left[-1,1\right]^{M}\backslash\left[-\delta,\delta\right]^{M}\right)<C\exp\left(-Dn\right),\quad n\in\IN.
\]

\section{\label{sec:Proofs}Proofs}

\subsection{Proof of Lemma \ref{lem:weak_conv}}

We show the implications $(i)\implies(ii)$, $(ii)\implies(iii)$,
and $(iii)\implies(i)$.

\subsubsection{$(i)\protect\implies(ii)$}

Let $\delta>0$ and $f:\IR^{M}\rightarrow\IR$ be a continuous and
bounded function with the properties
\[
f\geq0,\quad f(0)=0,\quad\textup{and}\quad f(x)=1,x\in\IR^{M}\backslash\prod_{\lambda=1}^{M}\left[-\delta,\delta\right].
\]
The weak convergence of $\left(\mu_{n}\right)_{n}$ to $\delta_{0}$
implies
\[
\int_{\IR^{M}}f\,\textup{d}\mu_{n}\xrightarrow[n\rightarrow\infty]{}\int_{\IR^{M}}f\,\textup{d}\delta_{0}=f(0)=0.
\]
Since $f\geq0$ and $f$ equals $1$ on $\IR^{M}\backslash\prod_{\lambda=1}^{M}\left[-\delta,\delta\right]$,
we have
\[
\mu_{n}\left(\IR^{M}\backslash\prod_{\lambda=1}^{M}\left[-\delta,\delta\right]\right)\xrightarrow[n\rightarrow\infty]{}0.
\]

\subsubsection{$(ii)\protect\implies(iii)$}

By $(ii)$, we have for all $m\in\IN$
\begin{equation}
\mu_{n}\left(\IR^{M}\backslash\prod_{\lambda=1}^{M}\left[-\frac{1}{m},\frac{1}{m}\right]\right)\xrightarrow[n\rightarrow\infty]{}0.\label{eq:mu_n_m}
\end{equation}
We define the natural numbers $n_{1},n_{2},\ldots$ as follows:
\[
n_{m}:=\min\left\{ k\in\IN\,\left|\,\mu_{n}\left(\prod_{\lambda=1}^{M}\left[-\frac{1}{m},\frac{1}{m}\right]\right)\geq1-\frac{1}{m},n\geq k\right.\right\} ,\quad m\in\IN.
\]
Note that the sets in the display are non-empty for each $m$ due
to (\ref{eq:mu_n_m}). Set
\[
\varepsilon_{n,\lambda}:=\frac{1}{m},\quad\lambda=1,\ldots,M,\quad n_{m}\leq n<n_{m+1},\quad m\in\IN.
\]
By construction, $\left\Vert \varepsilon_{n}\right\Vert _{\infty}\xrightarrow[n\rightarrow\infty]{}0$
and $\mu_{n}\left(\left[-\varepsilon_{n},\varepsilon_{n}\right]\right)\xrightarrow[n\rightarrow\infty]{}1$.

\subsubsection{$(iii)\protect\implies(i)$}

Let $f:\IR^{M}\rightarrow\IR$ be a continuous and bounded function.
Let $\eta>0$. Since $f$ is continuous, there is a $\delta>0$ such
that, for all $x\in\IR^{M}$,
\[
\left\Vert x\right\Vert _{\infty}<\delta\quad\textup{implies}\quad\left|f(x)-f(0)\right|<\eta.
\]
Then
\begin{equation}
\left|\int_{\IR^{M}}f\,\textup{d}\mu_{n}-\int_{\IR^{M}}f\,\textup{d}\delta_{0}\right|\leq\int_{\IR^{M}}\left|f(x)-f(0)\right|\,\mu_{n}\left(\textup{d}x\right).\label{eq:lem_disp_1}
\end{equation}
Let $n_{0}\in\IN$ be large enough that
\[
\left\Vert \varepsilon_{n}\right\Vert _{\infty}<\delta\quad\textup{and}\quad\mu_{n}\left(\left[-\varepsilon_{n},\varepsilon_{n}\right]\right)>1-\eta,\quad n\geq n_{0}.
\]
We continue with our calculation:
\begin{align*}
(\ref{eq:lem_disp_1}) & =\int_{\left[-\varepsilon_{n},\varepsilon_{n}\right]}\left|f(x)-f(0)\right|\,\mu_{n}\left(\textup{d}x\right)+\int_{\IR^{M}\backslash\left[-\varepsilon_{n},\varepsilon_{n}\right]}\left|f(x)-f(0)\right|\,\mu_{n}\left(\textup{d}x\right)\\
 & \leq\eta\mu_{n}\left(\left[-\varepsilon_{n},\varepsilon_{n}\right]\right)+2\left\Vert f\right\Vert _{\infty}\eta,
\end{align*}
which proves the weak convergence of $\left(\mu_{n}\right)_{n}$ to
$\delta_{0}$.

\subsection{Proof of Theorem \ref{thm:fast_conv}}

The proofs of the theorems in Section \ref{sec:Results} use characteristic
functions. Thus, we have to show pointwise convergence of the sequence
of characteristic functions of the normalised sums to the characteristic
function of the limiting distribution. In the case of Theorem \ref{thm:fast_conv},
we define $\varphi_{n}$ to be the characteristic function of the
distribution of
\[
\boldsymbol{S}_{n}=\left(\frac{S_{n,1}}{\sqrt{n_{1}}},\ldots,\frac{S_{n,M}}{\sqrt{n_{M}}}\right)
\]
and $\varphi_{\mathcal{N}\left(0,I_{M}\right)}$ to be the characteristic
function of $\mathcal{N}\left(0,I_{M}\right)$. Then our task is to
show $\varphi_{n}(t)\xrightarrow[n\rightarrow\infty]{}\varphi_{\mathcal{N}\left(0,I_{M}\right)}(t)$
for all $t\in\IR^{M}$. We will only prove the result for sequences
of de Finetti measures on $\left[-1,1\right]^{M}$. The results for
de Finetti sequences on $\IR^{M}$ can be shown analogously.

We will use the following notation for `asymptotic equivalence':

\begin{notation}Complex-valued sequences $\left(f_{n}\right)_{n},\left(g_{n}\right)_{n}$
are called \emph{asymptotically equal} (as $n\to\infty$), in short
$f_{n}\approx g_{n}$, if
\[
\lim_{n\rightarrow\infty}\frac{f_{n}}{g_{n}}=1.
\]
\end{notation}

Recall Definition \ref{def:Rademacher} and Notation \ref{not:prod_measure}.
By Definition \ref{def:de_Finetti_model}, $\varphi_{n}(t)$ can be
expressed as
\begin{align}
\IE\exp\left(it\cdot\boldsymbol{S}_{n}\right) & =\IE\exp\left(i\left(t_{1}\frac{S_{n,1}}{\sqrt{n_{1}}}+\cdots+t_{M}\frac{S_{n,M}}{\sqrt{n_{M}}}\right)\right)\nonumber \\
 & =\int_{\left[-1,1\right]^{M}}E_{m}\exp\left(i\left(t_{1}\frac{S_{n,1}}{\sqrt{n_{1}}}+\cdots+t_{M}\frac{S_{n,M}}{\sqrt{n_{M}}}\right)\right)\mu_{n}\left(\textup{d}m\right).\label{eq:1}
\end{align}

We first note that due to the boundedness of the integrand in (\ref{eq:1})
and the concentration property\\
$\mu_{n}\left[-\varepsilon_{n},\varepsilon_{n}\right]\xrightarrow[n\rightarrow\infty]{}1$,
we obtain
\begin{equation}
(\ref{eq:1})\approx\int_{\left[-\varepsilon_{n},\varepsilon_{n}\right]}E_{m}\exp\left(i\left(t_{1}\frac{S_{n,1}}{\sqrt{n_{1}}}+\cdots+t_{M}\frac{S_{n,M}}{\sqrt{n_{M}}}\right)\right)\mu_{n}\left(\textup{d}m\right).\label{eq:2}
\end{equation}
This step follows from the next lemma, provided that (\ref{eq:2})
converges to a non-zero limit, which we will show in (\ref{eq:2a}).
\begin{lem}
Let $f_{n}:\left[-1,1\right]^{M}\rightarrow\IC$ be measurable functions
with $\left\Vert f_{n}\right\Vert _{\infty}\leq1$, $n\in\IN$, and
$\left(\mu_{n}\right)_{n}$ the de Finetti sequence in the statement
of Theorem \ref{thm:fast_conv}. Then
\[
\left|\int_{\left[-1,1\right]^{M}}f_{n}\,\textup{d}\mu_{n}-\int_{\left[-\varepsilon_{n},\varepsilon_{n}\right]}f_{n}\,\textup{d}\mu_{n}\right|\xrightarrow[n\rightarrow\infty]{}0.
\]
\end{lem}

\begin{proof}
We estimate
\begin{align*}
\left|\int_{\left[-1,1\right]^{M}}f_{n}\,\textup{d}\mu_{n}-\int_{\left[-\varepsilon_{n},\varepsilon_{n}\right]}f_{n}\,\textup{d}\mu_{n}\right| & =\left|\int_{\IR^{M}\backslash\left[-\varepsilon_{n},\varepsilon_{n}\right]}f_{n}\,\textup{d}\mu_{n}\right|\leq\int_{\IR^{M}\backslash\left[-\varepsilon_{n},\varepsilon_{n}\right]}\left|f_{n}\right|\,\textup{d}\mu_{n}\\
 & \leq\mu_{n}\left(\IR^{M}\backslash\left[-\varepsilon_{n},\varepsilon_{n}\right]\right)\xrightarrow[n\rightarrow\infty]{}0.
\end{align*}
\end{proof}
Next, we have to calculate the expectation $E_{m}\exp\left(i\left(t_{1}S_{n,1}/\sqrt{n_{1}}+\cdots+t_{M}S_{n,M}/\sqrt{n_{M}}\right)\right)$.
Under $P_{m}^{\otimes n}$, all random variables $X_{\lambda i}$,
$\lambda=1,\ldots,M$, $i=1,\ldots,n_{\lambda}$ are independent.
Within each group $\lambda$, the $X_{\lambda i}$ are even i.i.d.
We subtract the expected value $m_{\lambda}$ of each $X_{\lambda i}$
under $P_{m}^{\otimes n}$. Note that for $t_{\lambda}=0$,
\[
E_{m}\exp\left(i\left(t_{\lambda}\frac{\sum_{j}\left(X_{\lambda j}-m_{\lambda}\right)}{\sqrt{n_{\lambda}}}\right)\right)=1.
\]
Now assume $t_{\lambda}\neq0$. Using a Taylor expansion yields
\begin{align*}
E_{m}\exp\left(i\left(t_{\lambda}\frac{\sum_{j}\left(X_{\lambda j}-m_{\lambda}\right)}{\sqrt{n_{\lambda}}}\right)\right) & =\left(1-\left(1-m_{\lambda}^{2}\right)\frac{t_{\lambda}^{2}}{2n_{\lambda}}+O\left(\frac{1}{n_{\lambda}^{3/2}}\right)\right)^{n_{\lambda}}.
\end{align*}
Thus, the conditional expectation in the integral in (\ref{eq:2})
can be written as
\[
\prod_{\lambda:\,t_{\lambda}\neq0}\exp\left(im_{\lambda}t_{\lambda}\sqrt{n_{\lambda}}\right)\left(1-\left(1-m_{\lambda}^{2}\right)\frac{t_{\lambda}^{2}}{2n_{\lambda}}+O\left(\frac{1}{n_{\lambda}^{3/2}}\right)\right)^{n_{\lambda}}.
\]
Set
\begin{align*}
a_{n,\lambda} & :=\exp\left(im_{\lambda}t_{\lambda}\sqrt{n_{\lambda}}\right),\\
b_{n,\lambda} & :=\left(1-\left(1-m_{\lambda}^{2}\right)\frac{t_{\lambda}^{2}}{2n_{\lambda}}+O\left(\frac{1}{n_{\lambda}^{3/2}}\right)\right)^{n_{\lambda}},\\
c_{n,\lambda} & :=a_{n,\lambda}b_{n,\lambda}
\end{align*}
for each $n$ and each $\lambda$.

First, we analyse the asymptotic behaviour of $\left(a_{n,\lambda}\right)_{n}$.
\begin{lem}
\label{lem:A_n}We have for all $n\in\IN$ and all $\lambda=1,\ldots,M$
\begin{align*}
a_{n,\lambda} & \in A_{n}:=\left\{ \exp\left(i\omega\right)\,\big|\,\omega\in\left[-\left|t_{\lambda}\right|\varepsilon_{n,\lambda}\sqrt{n_{\lambda}},\left|t_{\lambda}\right|\varepsilon_{n,\lambda}\sqrt{n_{\lambda}}\right]\right\} ,
\end{align*}
and
\[
A_{n}\searrow\left\{ 1\right\} \quad\textup{as }n\rightarrow\infty.
\]
\end{lem}

\begin{proof}
Since $t_{\lambda}\neq0$ and $m$ takes values on the interval $\left[-\varepsilon_{n},\varepsilon_{n}\right]$,
$a_{n,\lambda}$ lies on the unit circle in the complex plane, specifically
in the section defined by $A_{n}$. The convergence speed $\varepsilon_{n,\lambda}=o\left(\frac{1}{\sqrt{n_{\lambda}}}\right)$
implies $\pm\left|t_{\lambda}\right|\varepsilon_{n,\lambda}\sqrt{n_{\lambda}}\xrightarrow[n\rightarrow\infty]{}0$
and therefore $A_{n}\searrow\left\{ 1\right\} $ as $n\rightarrow\infty$.
\end{proof}
Next, we deal with the sequence $\left(b_{n,\lambda}\right)_{n}$.
We estimate
\begin{equation}
\left|b_{n,\lambda}-\exp\left(-\frac{t_{\lambda}^{2}}{2}\right)\right|\leq\left|b_{n,\lambda}-\left(1-\left(1-m_{\lambda}^{2}\right)\frac{t_{\lambda}^{2}}{2n_{\lambda}}\right)^{n_{\lambda}}\right|+\left|\left(1-\left(1-m_{\lambda}^{2}\right)\frac{t_{\lambda}^{2}}{2n_{\lambda}}\right)^{n_{\lambda}}-\exp\left(-\left(1-m_{\lambda}^{2}\right)\frac{t_{\lambda}^{2}}{2}\right)\right|,\label{eq:UB}
\end{equation}
and we will find upper bounds for each of the summands in (\ref{eq:UB}).
\begin{lem}
\label{lem:UB_1}There is a constant $C>0$ such that for all $n\in\IN$
we have
\[
\left|b_{n,\lambda}-\left(1-\left(1-m_{\lambda}^{2}\right)\frac{t_{\lambda}^{2}}{2n_{\lambda}}\right)^{n_{\lambda}}\right|\leq\frac{C}{n_{\lambda}^{3/2}}.
\]
\end{lem}

\begin{proof}
We define $r_{n_{\lambda}}:=\left(1-m_{\lambda}^{2}\right)\frac{t_{\lambda}^{2}}{2}$
and $s_{n_{\lambda}}$ to be the $O\left(\frac{1}{n_{\lambda}^{3/2}}\right)$
term in the definition of $b_{n,\lambda}$. Then we have
\[
b_{n,\lambda}=\left(1-\frac{r_{n_{\lambda}}}{n_{\lambda}}+s_{n_{\lambda}}\right)^{n_{\lambda}}=\sum_{k=0}^{n_{\lambda}}\left(1-\frac{r_{n_{\lambda}}}{n_{\lambda}}\right)^{n_{\lambda}-k}s_{n_{\lambda}}^{k}=\left(1-\frac{r_{n_{\lambda}}}{n_{\lambda}}\right)^{n_{\lambda}}\sum_{k=0}^{n_{\lambda}}\left(\frac{s_{n_{\lambda}}}{1-\frac{r_{n_{\lambda}}}{n_{\lambda}}}\right)^{k}.
\]
We calculate
\begin{align*}
\left|b_{n,\lambda}-\left(1-\left(1-m_{\lambda}^{2}\right)\frac{t_{\lambda}^{2}}{2n_{\lambda}}\right)^{n_{\lambda}}\right| & =\left|\left(1-\frac{r_{n_{\lambda}}}{n_{\lambda}}\right)^{n_{\lambda}}\sum_{k=0}^{n_{\lambda}}\left(\frac{s_{n_{\lambda}}}{1-\frac{r_{n_{\lambda}}}{n_{\lambda}}}\right)^{k}-\left(1-\frac{r_{n_{\lambda}}}{n_{\lambda}}\right)^{n_{\lambda}}\right|\\
 & =\left|1-\frac{r_{n_{\lambda}}}{n_{\lambda}}\right|^{n_{\lambda}}\left|\sum_{k=0}^{n_{\lambda}}\left(\frac{s_{n_{\lambda}}}{1-\frac{r_{n_{\lambda}}}{n_{\lambda}}}\right)^{k}-1\right|\\
 & =\left|1-\frac{r_{n_{\lambda}}}{n_{\lambda}}\right|^{n_{\lambda}}\left|\frac{1-\left(\frac{s_{n_{\lambda}}}{1-\frac{r_{n_{\lambda}}}{n_{\lambda}}}\right)^{n_{\lambda}+1}}{1-\frac{s_{n_{\lambda}}}{1-\frac{r_{n_{\lambda}}}{n_{\lambda}}}}-1\right|,
\end{align*}
where in the last step we applied the formula for the value of a geometric
series with complex summands. We note that $r_{n_{\lambda}}\xrightarrow[n_{\lambda}\rightarrow\infty]{}\frac{t_{\lambda}^{2}}{2}$
and hence
\[
\frac{s_{n_{\lambda}}}{1-\frac{r_{n_{\lambda}}}{n_{\lambda}}}\xrightarrow[n_{\lambda}\rightarrow\infty]{}0\quad\textup{and}\quad\left|1-\frac{r_{n_{\lambda}}}{n_{\lambda}}\right|^{n_{\lambda}}\xrightarrow[n_{\lambda}\rightarrow\infty]{}\exp\left(-\frac{t_{\lambda}^{2}}{2}\right).
\]
There is a constant $C_{1}>0$ such that
\[
\left|1-\frac{r_{n_{\lambda}}}{n_{\lambda}}\right|^{n_{\lambda}}\leq C_{1},\quad n_{\lambda}\in\IN.
\]
We continue with
\begin{align*}
\left|\frac{1-\left(\frac{s_{n_{\lambda}}}{1-\frac{r_{n_{\lambda}}}{n_{\lambda}}}\right)^{n_{\lambda}+1}}{1-\frac{s_{n_{\lambda}}}{1-\frac{r_{n_{\lambda}}}{n_{\lambda}}}}-1\right| & =\left|\frac{1-\left(\frac{s_{n_{\lambda}}}{1-\frac{r_{n_{\lambda}}}{n_{\lambda}}}\right)^{n_{\lambda}+1}-1+\frac{s_{n_{\lambda}}}{1-\frac{r_{n_{\lambda}}}{n_{\lambda}}}}{1-\frac{s_{n_{\lambda}}}{1-\frac{r_{n_{\lambda}}}{n_{\lambda}}}}\right|\\
 & =\left|\frac{s_{n_{\lambda}}}{1-\frac{r_{n_{\lambda}}}{n_{\lambda}}}\right|\left|\frac{1-\left(\frac{s_{n_{\lambda}}}{1-\frac{r_{n_{\lambda}}}{n_{\lambda}}}\right)^{n_{\lambda}}}{1-\frac{s_{n_{\lambda}}}{1-\frac{r_{n_{\lambda}}}{n_{\lambda}}}}\right|.
\end{align*}

Now we note that
\[
1-\frac{r_{n_{\lambda}}}{n_{\lambda}}\xrightarrow[n_{\lambda}\rightarrow\infty]{}1,\quad1-\frac{s_{n_{\lambda}}}{1-\frac{r_{n_{\lambda}}}{n_{\lambda}}}\xrightarrow[n_{\lambda}\rightarrow\infty]{}1,\quad\textup{and}\quad1-\left(\frac{s_{n_{\lambda}}}{1-\frac{r_{n_{\lambda}}}{n_{\lambda}}}\right)^{n_{\lambda}}\xrightarrow[n_{\lambda}\rightarrow\infty]{}1.
\]
The above implies
\[
1-\frac{r_{n_{\lambda}}}{n_{\lambda}},1-\frac{s_{n_{\lambda}}}{1-\frac{r_{n_{\lambda}}}{n_{\lambda}}},1-\left(\frac{s_{n_{\lambda}}}{1-\frac{r_{n_{\lambda}}}{n_{\lambda}}}\right)^{n_{\lambda}}=\Theta(1)
\]
and the existence of a constant $C_{2}>0$ such that
\[
\left|\frac{s_{n_{\lambda}}}{1-\frac{r_{n_{\lambda}}}{n_{\lambda}}}\right|\left|\frac{1-\left(\frac{s_{n_{\lambda}}}{1-\frac{r_{n_{\lambda}}}{n_{\lambda}}}\right)^{n_{\lambda}}}{1-\frac{s_{n_{\lambda}}}{1-\frac{r_{n_{\lambda}}}{n_{\lambda}}}}\right|\leq C_{2}\left|\frac{s_{n_{\lambda}}}{1-\frac{r_{n_{\lambda}}}{n_{\lambda}}}\right|,\quad n_{\lambda}\in\IN.
\]
Since $1-\frac{r_{n_{\lambda}}}{n_{\lambda}}\xrightarrow[n_{\lambda}\rightarrow\infty]{}1$,
there is a constant $C_{3}>0$ such that
\[
\left|\frac{s_{n_{\lambda}}}{1-\frac{r_{n_{\lambda}}}{n_{\lambda}}}\right|\leq C_{3}\left|s_{n_{\lambda}}\right|,\quad n_{\lambda}\in\IN.
\]
Finally, since $\left|s_{n_{\lambda}}\right|=O\left(\frac{1}{n_{\lambda}^{3/2}}\right)$,
there is a constant $C_{4}>0$ such that
\[
\left|s_{n_{\lambda}}\right|\leq\frac{C_{4}}{n_{\lambda}^{3/2}},\quad n_{\lambda}\in\IN.
\]
The claim follows by setting $C:=C_{1}C_{2}C_{3}C_{4}>0$.
\end{proof}
Next, we will find an upper bound for the second summand in (\ref{eq:UB}).
\begin{lem}
\label{lem:UB_2}There is a constant $D>0$ such that for all $n\in\IN$
we have
\[
\left|\left(1-\left(1-m_{\lambda}^{2}\right)\frac{t_{\lambda}^{2}}{2n_{\lambda}}\right)^{n_{\lambda}}-\exp\left(-\left(1-m_{\lambda}^{2}\right)\frac{t_{\lambda}^{2}}{2}\right)\right|\leq\frac{D}{n_{\lambda}}.
\]
\end{lem}

\begin{proof}
We again define $r_{n_{\lambda}}:=\left(1-m_{\lambda}^{2}\right)\frac{t_{\lambda}^{2}}{2}$.
Recall that $r_{n_{\lambda}}\xrightarrow[n_{\lambda}\rightarrow\infty]{}\frac{t_{\lambda}^{2}}{2}$.
Choose $m_{0}\in\IN$ large enough that
\begin{equation}
\left|\frac{r_{n_{\lambda}}}{n_{\lambda}}\right|<\frac{1}{2}\quad\textup{and}\quad\left|r_{n_{\lambda}}-\frac{t_{\lambda}^{2}}{2}\right|+2\frac{r_{n_{\lambda}}^{2}}{n_{\lambda}}<1,\quad n_{\lambda}\geq m_{0}.\label{eq:remainder_conds}
\end{equation}
We use Taylor expansions for the function $\exp$ and $\ln$:
\begin{align*}
\exp z & =1+z+O\left(z^{2}\right),\quad z\in\IR,\left|z\right|\leq1,\\
\ln\left(1+z\right) & =z+O\left(z^{2}\right),\quad z\in\IR,\left|z\right|\leq\frac{1}{2},
\end{align*}
More precisely, we have the bounds
\begin{align}
\left|\exp z-\left(1+z\right)\right| & \leq\frac{e}{2}z^{2},\quad z\in\IR,\left|z\right|\leq1,\label{eq:remainder_exp}\\
\left|\ln\left(1+z\right)-z\right| & \leq2z^{2},\quad z\in\IR,\left|z\right|\leq\frac{1}{2},\label{eq:remainder_ln}
\end{align}
on the remainder terms of order 2, as is easy to see from the Lagrange
form of said remainder terms. The constant $m_{0}$ in (\ref{eq:remainder_conds})
is chosen such that the value of $z$ in the Taylor expansions will
be small enough in our calculations.

Then we have for all $n_{\lambda}\geq m_{0}$
\begin{align*}
\left(1-\left(1-m_{\lambda}^{2}\right)\frac{t_{\lambda}^{2}}{2n_{\lambda}}\right)^{n_{\lambda}} & =\left(1-\frac{r_{n_{\lambda}}}{n_{\lambda}}\right)^{n_{\lambda}}=\exp\left[n_{\lambda}\ln\left(1-\frac{r_{n_{\lambda}}}{n_{\lambda}}\right)\right]\\
 & =\exp\left[n_{\lambda}\left(-\frac{r_{n_{\lambda}}}{n_{\lambda}}+s_{n_{\lambda}}\right)\right]=\exp\left[-r_{n_{\lambda}}+n_{\lambda}s_{n_{\lambda}}\right],
\end{align*}
where the remainder term $s_{n_{\lambda}}$ satisfies $\left|s_{n_{\lambda}}\right|\leq2\left(\frac{r_{n_{\lambda}}}{n_{\lambda}}\right)^{2}=O\left(\frac{1}{n_{\lambda}^{2}}\right)$
by (\ref{eq:remainder_ln}). We continue with our calculation:
\begin{align*}
\exp\left(-r_{n_{\lambda}}+n_{\lambda}s_{n_{\lambda}}\right) & =\exp\left[-\frac{t_{\lambda}^{2}}{2}-\left(r_{n_{\lambda}}-\frac{t_{\lambda}^{2}}{2}\right)+n_{\lambda}s_{n_{\lambda}}\right]\\
 & =\exp\left(-\frac{t_{\lambda}^{2}}{2}\right)\exp\left[-\left(r_{n_{\lambda}}-\frac{t_{\lambda}^{2}}{2}\right)+n_{\lambda}s_{n_{\lambda}}\right]\\
 & =\exp\left(-\frac{t_{\lambda}^{2}}{2}\right)\left[1-\left(r_{n_{\lambda}}-\frac{t_{\lambda}^{2}}{2}\right)+n_{\lambda}s_{n_{\lambda}}+u_{n_{\lambda}}\right],
\end{align*}
where the remainder term $u_{n_{\lambda}}$ satisfies
\begin{align}
\left|u_{n_{\lambda}}\right| & \leq\frac{e}{2}\left(\left(-\left(r_{n_{\lambda}}-\frac{t_{\lambda}^{2}}{2}\right)+n_{\lambda}s_{n_{\lambda}}\right)^{2}\right)\nonumber \\
 & =O\left(\left(-\left(r_{n_{\lambda}}-\frac{t_{\lambda}^{2}}{2}\right)+n_{\lambda}s_{n_{\lambda}}\right)^{2}\right)\nonumber \\
 & =O\left(\left(r_{n_{\lambda}}-\frac{t_{\lambda}^{2}}{2}\right)^{2}\right)+O\left(\left(n_{\lambda}s_{n_{\lambda}}\right)^{2}\right)=O\left(\left(r_{n_{\lambda}}-\frac{t_{\lambda}^{2}}{2}\right)^{2}\right)+O\left(\frac{1}{n_{\lambda}^{2}}\right),\label{eq:u_n}
\end{align}
 by (\ref{eq:remainder_exp}). We used $O\left(\left(n_{\lambda}s_{n_{\lambda}}\right)^{2}\right)=O\left(\frac{n_{\lambda}^{2}}{n_{\lambda}^{4}}\right)=O\left(\frac{1}{n_{\lambda}^{2}}\right)$
in the last step above. Now we estimate the difference
\begin{align*}
\left|\left(1-\left(1-m_{\lambda}^{2}\right)\frac{t_{\lambda}^{2}}{2n_{\lambda}}\right)^{n_{\lambda}}-\exp\left(-\left(1-m_{\lambda}^{2}\right)\frac{t_{\lambda}^{2}}{2}\right)\right| & =\left|\left(1-\frac{r_{n_{\lambda}}}{n_{\lambda}}\right)^{n_{\lambda}}-\exp\left(-\frac{t_{\lambda}^{2}}{2}\right)\right|\\
 & =\exp\left(-\frac{t_{\lambda}^{2}}{2}\right)\left|1-\left(r_{n_{\lambda}}-\frac{t_{\lambda}^{2}}{2}\right)+n_{\lambda}s_{n_{\lambda}}+u_{n_{\lambda}}-1\right|\\
 & =\exp\left(-\frac{t_{\lambda}^{2}}{2}\right)\left|-\left(r_{n_{\lambda}}-\frac{t_{\lambda}^{2}}{2}\right)+n_{\lambda}s_{n_{\lambda}}+u_{n_{\lambda}}\right|
\end{align*}
Taking into account $n_{\lambda}s_{n_{\lambda}}=O\left(\frac{1}{n_{\lambda}}\right)$
and (\ref{eq:u_n}), we have
\begin{align}
\exp\left(-\frac{t_{\lambda}^{2}}{2}\right)\left|-\left(r_{n_{\lambda}}-\frac{t_{\lambda}^{2}}{2}\right)+n_{\lambda}s_{n_{\lambda}}+u_{n_{\lambda}}\right| & \leq\exp\left(-\frac{t_{\lambda}^{2}}{2}\right)\left[\left|r_{n_{\lambda}}-\frac{t_{\lambda}^{2}}{2}\right|+\left|n_{\lambda}s_{n_{\lambda}}\right|+\left|u_{n_{\lambda}}\right|\right]\nonumber \\
 & =O\left(r_{n_{\lambda}}-\frac{t_{\lambda}^{2}}{2}\right)+O\left(\frac{1}{n_{\lambda}}\right)+O\left(\left(r_{n_{\lambda}}-\frac{t_{\lambda}^{2}}{2}\right)^{2}\right)+O\left(\frac{1}{n_{\lambda}^{2}}\right)\nonumber \\
 & =O\left(r_{n_{\lambda}}-\frac{t_{\lambda}^{2}}{2}\right)+O\left(\frac{1}{n_{\lambda}}\right).\label{eq:diff}
\end{align}
Using the definition of $r_{n_{\lambda}}$ and the assumption $\varepsilon_{n,\lambda}=o\left(\frac{1}{\sqrt{n_{\lambda}}}\right)$,
\begin{align*}
\left|r_{n_{\lambda}}-\frac{t_{\lambda}^{2}}{2}\right| & =m_{\lambda}^{2}\frac{t_{\lambda}^{2}}{2}\leq\varepsilon_{n,\lambda}^{2}\frac{t_{\lambda}^{2}}{2}=o\left(\frac{1}{n_{\lambda}}\right).
\end{align*}
It follows from the above and (\ref{eq:diff}) that
\[
\left|\left(1-\left(1-m_{\lambda}^{2}\right)\frac{t_{\lambda}^{2}}{2n_{\lambda}}\right)^{n_{\lambda}}-\exp\left(-\left(1-m_{\lambda}^{2}\right)\frac{t_{\lambda}^{2}}{2}\right)\right|=O\left(\frac{1}{n_{\lambda}}\right)
\]
and hence the statement of the lemma is proved.
\end{proof}
Lemmas \ref{lem:UB_1} and \ref{lem:UB_2} imply the following corollary:
\begin{cor}
\label{cor:B_n}The range of $b_{n,\lambda}$ over all possible values
of $m_{\lambda}\in\left[-\varepsilon_{n,\lambda},\varepsilon_{n,\lambda}\right]$
lies in a $\frac{D}{n_{\lambda}}$-neighbourhood $B_{n}$ of the set
$\left[\exp\left(-\frac{t_{\lambda}^{2}}{2}\right),\exp\left(-\left(1-\varepsilon_{n,\lambda}^{2}\right)\frac{t_{\lambda}^{2}}{2}\right)\right]$
in $\IC$. The set sequence $B_{n}$ converges:
\[
B_{n}\searrow\left\{ \exp\left(-\frac{t_{\lambda}^{2}}{2}\right)\right\} \quad\textup{as }n\rightarrow\infty.
\]
\end{cor}

Since $\mu_{n}$ is symmetric, once we integrate $c_{n,\lambda}$
with respect to $\mu_{n}$, only the integral of the real part $\textup{Re}\left(c_{n,\lambda}\right)$
remains. This real part is the radius of the circle on which all values
of $c_{n,\lambda}$ lie, and these radii are the values that comprise
the set $B_{n}$. As each $\mu_{n}$ is a probability measure and
we assumed $\mu_{n}\left[-\varepsilon_{n},\varepsilon_{n}\right]\xrightarrow[n\rightarrow\infty]{}1$,
we have
\begin{align}
(\ref{eq:2}) & \xrightarrow[n\rightarrow\infty]{}\prod_{\lambda=1}^{M}\exp\left(-\frac{t_{\lambda}^{2}}{2}\right)=\varphi_{\mathcal{N}\left(0,I_{M}\right)}(t).\label{eq:2a}
\end{align}
As we see, (\ref{eq:2}) converges to a non-zero limit as $n\rightarrow\infty$,
thereby proving the asymptotic equivalence of (\ref{eq:2}) to $\varphi_{n}(t)$.

This concludes the proof of Theorem \ref{thm:fast_conv}.

\subsection{Proof of Theorem \ref{thm:phase_transition}}

We now show the statements of Theorem \ref{thm:phase_transition}.
The first statement is a corollary of Theorem \ref{thm:fast_conv}.
To prove the second statement, we once again have to calculate $\varphi_{n}(t)=\IE\exp\left(it\cdot\boldsymbol{S}_{n}\right)$.
The normalisation $\gamma_{n,\lambda}=\sqrt{n_{\lambda}}$ is the
same as in Theorem \ref{thm:fast_conv}. Since $\mu_{n}$ is the contraction
of $\mu$, which has support in $\left[-1,1\right]^{M}$, the support
of $\mu_{n}$ is a subset of $\left[-\varepsilon_{n},\varepsilon_{n}\right]$
for each $n\in\IN$, and hence
\begin{equation}
\varphi_{n}(t)=\int_{\left[-\varepsilon_{n},\varepsilon_{n}\right]}E_{m}\exp\left(i\left(t_{1}\frac{S_{n,1}}{\sqrt{n_{1}}}+\cdots+t_{M}\frac{S_{n,M}}{\sqrt{n_{M}}}\right)\right)\mu_{n}\left(\textup{d}m\right).\label{eq:2.1}
\end{equation}
At this point, we change variables by setting $s:=\left(\frac{1}{\varepsilon_{n,1}},\ldots,\frac{1}{\varepsilon_{n,M}}\right)\circ m$.
Then
\begin{align}
(\ref{eq:2.1}) & =\int_{\left[-1,1\right]^{M}}E_{\varepsilon_{n}\circ s}\exp\left(i\left(t_{1}\frac{S_{n,1}}{\sqrt{n_{1}}}+\cdots+t_{M}\frac{S_{n,M}}{\sqrt{n_{M}}}\right)\right)\mu_{n}\left(\textup{d}\left(\varepsilon_{n}\circ s\right)\right)\nonumber \\
 & =\int_{\left[-1,1\right]^{M}}E_{\varepsilon_{n}\circ s}\exp\left(i\left(t_{1}\frac{S_{n,1}}{\sqrt{n_{1}}}+\cdots+t_{M}\frac{S_{n,M}}{\sqrt{n_{M}}}\right)\right)\mu\left(\textup{d}s\right),\label{eq:2.2}
\end{align}
where we used that, by definition, $\mu_{n}\left(\varepsilon_{n}\circ A\right)=\mu\left(A\right)$
for all measurable sets $A$. The conditional expectations $E_{\varepsilon_{n}\circ s}\exp\left(it_{\lambda}\frac{S_{n,\lambda}}{\sqrt{n_{\lambda}}}\right)$
can be expressed via a Taylor expansion as
\begin{align*}
E_{\varepsilon_{n}\circ s}\exp\left(it_{\lambda}\frac{S_{n,\lambda}}{\sqrt{n_{\lambda}}}\right) & =\exp\left(i\varepsilon_{n,\lambda}s_{\lambda}t_{\lambda}\sqrt{n_{\lambda}}\right)\left(1-\left(1-\varepsilon_{n,\lambda}^{2}s_{\lambda}^{2}\right)\frac{t_{\lambda}^{2}}{2n_{\lambda}}+O\left(\frac{1}{n_{\lambda}^{3/2}}\right)\right)^{n_{\lambda}}\\
 & \xrightarrow[n\rightarrow\infty]{}\exp\left(it_{\lambda}h_{\lambda}s_{\lambda}\right)\exp\left(-\frac{t_{\lambda}^{2}}{2}\right).
\end{align*}
By dominated convergence, we obtain
\begin{align*}
(\ref{eq:2.2}) & \xrightarrow[n\rightarrow\infty]{}\int_{\left[-1,1\right]^{M}}\prod_{\lambda=1}^{M}\exp\left(it_{\lambda}h_{\lambda}s_{\lambda}\right)\exp\left(-\frac{t_{\lambda}^{2}}{2}\right)\mu\left(\textup{d}s\right)=\varphi_{\mu}\left(h\circ t\right)\varphi_{\mathcal{N}\left(0,I_{M}\right)}(t),
\end{align*}
where $\varphi_{\mu}$ is the characteristic function of $\mu$.

Since $\varphi_{\mu}\left(h\circ t\right)=\varphi_{h\circ\mu}\left(t\right)$,
this concludes the proof of the second statement of Theorem \ref{thm:phase_transition}.

In statement 3, we have the normalisation sequences $\gamma_{n,\lambda}:=\varepsilon_{n,\lambda}n_{\lambda}$.
Let $\varphi_{n}$ be the characteristic function of the distribution
of
\[
\boldsymbol{S}_{n}=\left(\frac{S_{n,1}}{\gamma_{n,1}},\ldots,\frac{S_{n,M}}{\gamma_{n,M}}\right).
\]
The proof of statement 3 proceeds as outlined above and yields
\begin{equation}
\varphi_{n}(t)=\int_{\left[-1,1\right]^{M}}E_{\varepsilon_{n}\circ s}\exp\left(i\left(t_{1}\frac{S_{n,1}}{\gamma_{n,1}}+\cdots+t_{M}\frac{S_{n,M}}{\gamma_{n,M}}\right)\right)\mu\left(\textup{d}s\right).\label{eq:3.1}
\end{equation}
We once again expand the conditional expectation in the integral above:
\begin{align*}
E_{\varepsilon_{n}\circ s}\exp\left(it_{\lambda}\frac{S_{n,\lambda}}{\gamma_{n,\lambda}}\right) & =\exp\left(i\varepsilon_{n,\lambda}s_{\lambda}t_{\lambda}\frac{n_{\lambda}}{\gamma_{n,\lambda}}\right)\left(1-\left(1-\varepsilon_{n,\lambda}^{2}s_{\lambda}^{2}\right)\frac{t_{\lambda}^{2}}{2\gamma_{n,\lambda}^{2}}+O\left(\frac{1}{\gamma_{n,\lambda}^{3}}\right)\right)^{n_{\lambda}}.
\end{align*}
Recall that we assumed $1/\sqrt{n_{\lambda}}=o\left(\varepsilon_{n,\lambda}\right)$.
It follows that $1/n_{\lambda}=o\left(\varepsilon_{n,\lambda}^{2}\right)$
and $n_{\lambda}\varepsilon_{n,\lambda}^{2}\xrightarrow[n\rightarrow\infty]{}\infty$.
Hence, $n_{\lambda}=o\left(\varepsilon_{n,\lambda}^{2}n_{\lambda}^{2}\right)$,
and due to $\gamma_{n,\lambda}^{2}=\varepsilon_{n,\lambda}^{2}n_{\lambda}^{2}$
we have
\[
E_{\varepsilon_{n}\circ s}\exp\left(it_{\lambda}\frac{S_{n,\lambda}}{\gamma_{n,\lambda}}\right)\approx\exp\left(i\varepsilon_{n,\lambda}s_{\lambda}t_{\lambda}\frac{n_{\lambda}}{\gamma_{n,\lambda}}\right)=\exp\left(is_{\lambda}t_{\lambda}\right).
\]

By dominated convergence,
\begin{align*}
\varphi_{n}(t) & \xrightarrow[n\rightarrow\infty]{}\int_{\left[-1,1\right]^{M}}\prod_{\lambda=1}^{M}\exp\left(is_{\lambda}t_{\lambda}\right)\mu\left(\textup{d}s\right)=\varphi_{\mu}\left(t\right),
\end{align*}
and our proof is complete.

We omit the proof of Theorem \ref{thm:clusters} as it amounts to
a repetition of the arguments presented above, carefully keeping track
of the coordinates belonging to each of the three clusters. Corollaries
\ref{cor:voting_margins} and \ref{cor:single_group} follow immediately
from Theorem \ref{thm:clusters}.

\subsection{Proof of Proposition \ref{prop:corr}}

We calculate the correlation $\IE\left(X_{\lambda1}X_{\nu2}\right)$
using Definition \ref{def:de_Finetti_model}:
\begin{align*}
\IE\left(X_{\lambda1}X_{\nu2}\right) & =\int_{\IR^{M}}E_{\bar{m}}\left(X_{\lambda1}X_{\nu2}\right)\,\mu_{n}\left(\textup{d}m\right)=\int_{\IR^{M}}\left(P_{\bar{m}}^{\otimes n}\left(X_{\lambda1}=X_{\nu2}\right)-P_{\bar{m}}^{\otimes n}\left(X_{\lambda1}\neq X_{\nu2}\right)\right)\,\mu_{n}\left(\textup{d}m\right)\\
 & =2\int_{\IR^{M}}P_{\bar{m}}^{\otimes n}\left(X_{\lambda1}=X_{\nu2}\right)\,\mu_{n}\left(\textup{d}m\right)-1.
\end{align*}
We express $P_{\bar{m}}^{\otimes n}\left(X_{\lambda1}=X_{\nu2}\right)$
in terms of $\bar{m}$:
\begin{align*}
P_{\bar{m}}^{\otimes n}\left(X_{\lambda1}=X_{\nu2}\right) & =P_{\bar{m}}^{\otimes n}\left(X_{\lambda1}=X_{\nu2}=1\right)+P_{\bar{m}}^{\otimes n}\left(X_{\lambda1}=X_{\nu2}=-1\right)\\
 & =\frac{1+\bar{m}_{\lambda}}{2}\frac{1+\bar{m}_{\nu}}{2}+\frac{1-\bar{m}_{\lambda}}{2}\frac{1-\bar{m}_{\nu}}{2}=\frac{1+\bar{m}_{\lambda}\bar{m}_{\nu}}{2},
\end{align*}
where in the second step we used that all $X_{\kappa i}$ are independent
random variables under $P_{\bar{m}}^{\otimes n}$. Hence,
\[
\IE\left(X_{\lambda1}X_{\nu2}\right)=\int_{\IR^{M}}\left(1+\bar{m}_{\lambda}\bar{m}_{\nu}\right)\,\mu_{n}\left(\textup{d}m\right)-1=\int_{\IR^{M}}\bar{m}_{\lambda}\bar{m}_{\nu}\,\mu_{n}\left(\textup{d}m\right).
\]
We have
\begin{align*}
\int_{\IR^{M}}\bar{m}_{\lambda}\bar{m}_{\nu}\,\mu_{n}\left(\textup{d}m\right) & =\int_{\left[-\varepsilon_{n},\varepsilon_{n}\right]}\bar{m}_{\lambda}\bar{m}_{\nu}\,\mu_{n}\left(\textup{d}m\right)+\int_{\IR^{M}\backslash\left[-\varepsilon_{n},\varepsilon_{n}\right]}\bar{m}_{\lambda}\bar{m}_{\nu}\,\mu_{n}\left(\textup{d}m\right),
\end{align*}
and we upper bound each summand in turn:
\[
\left|\int_{\left[-\varepsilon_{n},\varepsilon_{n}\right]}\bar{m}_{\lambda}\bar{m}_{\nu}\,\mu_{n}\left(\textup{d}m\right)\right|\leq\sup_{m\in\left[-\varepsilon_{n},\varepsilon_{n}\right]}\left|\bar{m}_{\lambda}\bar{m}_{\nu}\right|=O\left(\varepsilon_{n,\lambda}\varepsilon_{n,\nu}\right)\xrightarrow[n\rightarrow\infty]{}0,
\]
where we used the fact that $\lim_{m\rightarrow0}\frac{m_{\kappa}}{\bar{m}_{\kappa}}=1$
for all $\kappa=1,\ldots,M$, and
\[
\left|\int_{\IR^{M}\backslash\left[-\varepsilon_{n},\varepsilon_{n}\right]}\bar{m}_{\lambda}\bar{m}_{\nu}\,\mu_{n}\left(\textup{d}m\right)\right|\leq\mu_{n}\left(\IR^{M}\backslash\left[-\varepsilon_{n},\varepsilon_{n}\right]\right)\xrightarrow[n\rightarrow\infty]{}0,
\]
where we used that $\left\Vert \bar{m}_{\lambda}\bar{m}_{\nu}\right\Vert _{\infty}\leq1$.

\subsection{Proof of Theorem \ref{thm:LLT}}

For our proof of the local limit Theorem \ref{thm:LLT}, we will need
the following auxiliary lemmas:
\begin{lem}
\label{lem:char_fn}Let $Y:=\left(Y_{1},\ldots,Y_{d}\right)$ be a
random vector with values on the lattice $\prod_{\lambda=1}^{d}\left(v_{\lambda}+w_{\lambda}\mathbb{Z}\right)$,
where $v_{\lambda}\in\IR$ and $w_{\lambda}>0$, and let $\varphi$
be the characteristic function of the distribution of $Y$. Then we
have for all $k_{1},\ldots,k_{d}\in\mathbb{Z}$
\[
\left|\varphi\left(2\pi\left(\frac{k_{1}}{w_{1}},\cdots,\frac{k_{d}}{w_{d}}\right)\right)\right|=1,
\]
and for all $t\in\mathbb{R}^{d}$ such that $0<t_{\lambda}<\frac{2\pi}{w_{\lambda}}$
for at least one component $\lambda$, we have
\[
\left|\varphi\left(t\right)\right|<1.
\]
\end{lem}

\begin{proof}
The lemma follows from the proof of Theorem 3.5.2 on page 140 in \cite{Durrett}.
\end{proof}
Lemma \ref{lem:char_fn} gives an upper bound for the characteristic
function of a random vector on a lattice (such as $\boldsymbol{S}_{n}$),
which we shall use in our calculations later on. We will use the following
inversion formulas to recover distributions from their characteristic
functions:
\begin{lem}
\label{lem:discr_inv}Let $\left(Y_{1},\ldots,Y_{d}\right)$ be a
random vector as in Lemma \ref{lem:char_fn} with characteristic function
$\varphi$. Then, for all $x\in\prod_{\lambda=1}^{d}\left(v_{\lambda}+w_{\lambda}\mathbb{Z}\right)$,
\[
\mathbb{P}\left(\left(Y_{1},\ldots,Y_{d}\right)=x\right)=\frac{\prod_{\lambda=1}^{d}w_{\lambda}}{\left(2\pi\right)^{d}}\int_{\prod_{\lambda}\left[-\frac{\pi}{w_{\lambda}},\frac{\pi}{w_{\lambda}}\right]}e^{-it\cdot x}\,\varphi(t)\,\textup{d}t.
\]
\end{lem}

\begin{proof}
See e.g.\!  Section 3.10 in \cite{Durrett}.
\end{proof}
\begin{lem}
\label{lem:cont_inv}Let $\varphi$ be the characteristic function
of some probability distribution on $\IR^{d}$ such that $\varphi$
is Lebesgue integrable. Then
\[
f(x)=\frac{1}{\left(2\pi\right)^{d}}\int_{\mathbb{R}^{d}}e^{-it\cdot x}\,\varphi(t)\,\textup{d}t
\]
for all $x\in\IR^{d}$ defines a continuous Lebesgue density function
$f$ for said probability distribution.
\end{lem}

\begin{proof}
This is Theorem 5.5 in \cite{Wengenroth} (alternatively, see Theorem
3.3.14 in \cite{Durrett}).
\end{proof}
Now we start the proof proper of Theorem \ref{thm:LLT}. Let $\II_{A}$
be the indicator function of a measurable set $A$. The lattice on
which $\boldsymbol{S}_{n}$ lives is $\mathcal{L}_{n}:=\prod_{\lambda=1}^{M}\left(\sqrt{n_{\lambda}}+\frac{2}{\sqrt{n_{\lambda}}}\mathbb{Z}\right)$.
Let $x\in\mathcal{L}_{n}$. By Lemmas \ref{lem:discr_inv} and \ref{lem:cont_inv},
we have
\begin{align}
\left|\frac{\prod_{\lambda=1}^{M}\sqrt{n_{\lambda}}}{2^{M}}\mathbb{P}\left(\boldsymbol{S}_{n}=x\right)-\phi(x)\right| & =\frac{1}{\left(2\pi\right)^{M}}\left|\int_{\prod_{\lambda}\left[-\frac{\pi\sqrt{n_{\lambda}}}{2},\frac{\pi\sqrt{n_{\lambda}}}{2}\right]}e^{-it\cdot x}\varphi_{n}(t)\,\textup{d}t-\int_{\mathbb{R}^{M}}e^{-it\cdot x}\varphi_{\mathcal{N}\left(0,I_{M}\right)}(t)\,\textup{d}t\right|\nonumber \\
 & \leq\frac{1}{\left(2\pi\right)^{M}}\int_{\mathbb{R}^{M}}\II_{\prod_{\lambda}\left[-\frac{\pi\sqrt{n_{\lambda}}}{2},\frac{\pi\sqrt{n_{\lambda}}}{2}\right]}(t)\left|\varphi_{n}(t)-\varphi_{\mathcal{N}\left(0,I_{M}\right)}(t)\right|\textup{d}t\label{eq:centre}\\
 & \quad+\frac{1}{\left(2\pi\right)^{M}}\int_{\mathbb{R}^{M}\backslash\prod_{\lambda}\left[-\frac{\pi\sqrt{n_{\lambda}}}{2},\frac{\pi\sqrt{n_{\lambda}}}{2}\right]}\left|\varphi_{\mathcal{N}\left(0,I_{M}\right)}(t)\right|\textup{d}t.\label{eq:tails}
\end{align}
Note that both (\ref{eq:centre}) and (\ref{eq:tails}) are independent
of the point $x\in\mathcal{L}_{n}$. The term (\ref{eq:tails}) converges
to 0 as $n\rightarrow\infty$, since $\left|\varphi_{\mathcal{N}\left(0,I_{M}\right)}\right|$
is integrable. Thus, our remaining task is to show that (\ref{eq:centre})
converges to 0. By Theorem \ref{thm:fast_conv}, which we have already
shown, $\varphi_{n}(t)\xrightarrow[n\rightarrow\infty]{}\varphi_{\mathcal{N}\left(0,I_{M}\right)}(t)$
holds for all $t\in\IR^{M}$. Hence, we can prove that (\ref{eq:centre})
converges to 0 by finding an appropriate integrable majorant and applying
the theorem of dominated convergence.

We pick some $0<\delta<\pi/2$ and partition the set $\prod_{\lambda}\left[-\frac{\pi\sqrt{n_{\lambda}}}{2},\frac{\pi\sqrt{n_{\lambda}}}{2}\right]$
into the disjoint sets\\
$A_{n}:=\prod_{\lambda}\left[-\delta\sqrt{n_{\lambda}},\delta\sqrt{n_{\lambda}}\right]$
and $B_{n}:=\prod_{\lambda}\left[-\frac{\pi\sqrt{n_{\lambda}}}{2},\frac{\pi\sqrt{n_{\lambda}}}{2}\right]\backslash A_{n}$
for each $n\in\mathbb{N}$. Let for all $s\in\left[-1,1\right]$ $\varphi_{\mathcal{R}(s)}$
be the characteristic function of the Rademacher distribution on $\left\{ -1,1\right\} $
(see Definition \ref{def:Rademacher}).

We construct an integrable majorant first over $A_{n}$ and then $B_{n}$.
The following upper bound holds over $A_{n}$:
\begin{equation}
\II_{A_{n}}(t)\left|\varphi_{n}(t)\right|\leq\II_{A_{n}}(t)\int_{\left[-1,1\right]^{M}}\prod_{\lambda}\left|\varphi_{\mathcal{R}(m_{\lambda})}\left(\frac{t_{\lambda}}{\sqrt{n_{\lambda}}}\right)\right|^{n_{\lambda}}\mu_{n}(\textup{d}m).\label{eq:upper_bound_int}
\end{equation}
The next lemma will provide an upper bound for the characteristic
function of the Rademacher distribution.
\begin{lem}
\label{lem:UB_exponential}For any $x\in\IR$ and any $n\in\IN$,
the upper bound
\[
\left|\exp\left(ix\right)-\sum_{k=0}^{n}\frac{\left(ix\right)^{k}}{k!}\right|\leq\min\left\{ \frac{\left|x\right|^{n+1}}{\left(n+1\right)!},\frac{2\left|x\right|^{n}}{n!}\right\} 
\]
holds.
\end{lem}

\begin{proof}
We use Taylor expansions with the remainder term in Legendre form.
First we expand to order $n$,
\begin{equation}
\exp\left(ix\right)-\sum_{k=0}^{n}\frac{\left(ix\right)^{k}}{k!}=\frac{\left(ix\right)^{n+1}\exp\left(i\xi_{n+1}\right)}{\left(n+1\right)!},\label{eq:lem_UB_1}
\end{equation}
and then to order $n-1$,
\[
\exp\left(ix\right)-\sum_{k=0}^{n-1}\frac{\left(ix\right)^{k}}{k!}=\frac{\left(ix\right)^{n}\exp\left(i\xi_{n}\right)}{n!}.
\]
The variables $\xi_{n}$ and $\xi_{n+1}$ each lie between 0 and $x$.
Combining the two displays above yields
\begin{align}
\exp\left(ix\right)-\sum_{k=0}^{n}\frac{\left(ix\right)^{k}}{k!} & =\exp\left(ix\right)-\sum_{k=0}^{n-1}\frac{\left(ix\right)^{k}}{k!}-\frac{\left(ix\right)^{n}}{n!}\nonumber \\
 & =\frac{\left(ix\right)^{n}\exp\left(i\xi_{n}\right)}{n!}-\frac{\left(ix\right)^{n}}{n!}.\label{eq:lem_UB_2}
\end{align}
Taking into account $\left|i^{n}\right|=\left|i^{n+1}\right|=\left|\exp\left(i\xi_{n}\right)\right|=1$,
we obtain the desired result by applying upper bounds to the absolute
value of (\ref{eq:lem_UB_1}) and (\ref{eq:lem_UB_2}). 
\end{proof}
Applying the lemma to $x=X_{\lambda1}-s$ and $n=2$ and taking the
expectation $E_{s}$, we obtain
\begin{align}
\left|E_{s}\exp\left(iu\left(X_{\lambda1}-s\right)\right)-E_{s}\left(\sum_{k=0}^{2}\frac{\left(iu\right)^{k}\left(X_{\lambda1}-s\right)^{k}}{k!}\right)\right| & \leq u^{2}E_{s}\min\left\{ \frac{\left|u\right|\left|X_{\lambda1}-s\right|^{3}}{3!},\left|X_{\lambda1}-s\right|^{2}\right\} \nonumber \\
 & \leq u^{2}\min\left\{ \frac{\left|u\right|E_{s}\left|X_{\lambda1}-s\right|^{3}}{3!},E_{s}\left|X_{\lambda1}-s\right|^{2}\right\} .\label{eq:Rad_expec_UB}
\end{align}
We calculate the expectations in the upper bound above:
\[
E_{s}\left|X_{\lambda1}-s\right|^{2}=1-s^{2}\quad\textup{and}\quad E_{s}\left|X_{\lambda1}-s\right|^{3}=2\left|s\right|\left(1-s^{2}\right).
\]

We calculate an upper bound for the Rademacher characteristic function:
\begin{align*}
\left|\varphi_{\mathcal{R}(s)}(u)\right| & =\left|E_{s}\exp\left(iu\left(X_{\lambda1}-s\right)\right)\right|\left|\exp\left(ius\right)\right|\\
 & \leq\left|1-\left(1-s^{2}\right)\frac{u^{2}}{2}\right|+u^{2}\left(1-s^{2}\right)\min\left\{ \left|u\right|\left|s\right|,1\right\} \\
 & \leq1-\left(1-s^{2}\right)\frac{u^{2}}{2}+\left(1-s^{2}\right)\frac{u^{2}}{4}\leq\exp\left(-\left(1-s^{2}\right)\frac{u^{2}}{4}\right).
\end{align*}
Above, in the first inequality, we used a Taylor expansion of order
2 with remainder term upper bounded by (\ref{eq:Rad_expec_UB}). The
second inequality holds for all $u\in\IR$ such that $\left|u\right|\left|s\right|\leq1/4$,
for which we have
\[
\left|u\right|^{3}\left(1-s^{2}\right)\left|s\right|\leq\left(1-s^{2}\right)\frac{u^{2}}{4}.
\]
The third inequality follows from $1-x\leq\exp\left(-x\right)$ for
all $x\in\IR$.

Therefore,
\[
\left|\varphi_{\mathcal{R}\left(m_{\lambda}\right)}\left(\frac{t_{\lambda}}{\sqrt{n_{\lambda}}}\right)\right|\leq\exp\left(-\left(1-m_{\lambda}^{2}\right)\frac{t_{\lambda}^{2}}{4n_{\lambda}}\right)
\]
holds for all $n\in\IN$ large enough.

By assumption, there is some $\tau\in\left(0,1\right)^{M}$ with the
property (\ref{eq:summability}). We continue with our calculation:
\begin{align}
(\ref{eq:upper_bound_int}) & \leq\II_{A_{n}}(t)\int_{\left[-1,1\right]^{M}}\exp\left(-\frac{1}{4}\sum_{\lambda}\left(1-m_{\lambda}^{2}\right)t_{\lambda}^{2}\right)\mu_{n}(\textup{d}m)\nonumber \\
 & \leq\int_{\left[-\tau,\tau\right]}\exp\left(-\frac{1}{4}\sum_{\lambda}\left(1-m_{\lambda}^{2}\right)t_{\lambda}^{2}\right)\mu_{n}(\textup{d}m)+\II_{A_{n}}(t)\,\mu_{n}\left(\left[-1,1\right]^{M}\backslash\left[-\tau,\tau\right]\right)\nonumber \\
 & \leq\exp\left(-\frac{1}{4}\sum_{\lambda}\left(1-\tau_{\lambda}^{2}\right)t_{\lambda}^{2}\right)+\II_{A_{n}}(t)\,\mu_{n}\left(\left[-1,1\right]^{M}\backslash\left[-\tau,\tau\right]\right).\label{eq:upper_bound_int_2}
\end{align}
It is clear that the first summand in (\ref{eq:upper_bound_int_2})
is integrable. For the second summand, we have
\begin{align*}
\II_{A_{n}}(t)\,\mu_{n}\left(\left[-1,1\right]^{M}\backslash\left[-\tau,\tau\right]\right) & \leq\II_{A_{1}}\left(t\right)\mu_{1}\left(\left[-1,1\right]^{M}\backslash\left[-\tau,\tau\right]\right)+\sum_{k=1}^{\infty}\II_{A_{k+1}\backslash A_{k}}\left(t\right)\,\mu_{k}\left(\left[-1,1\right]^{M}\backslash\left[-\tau,\tau\right]\right)\\
 & =:f(t).
\end{align*}
Let $\boldsymbol{\lambda}^{M}$ be the Lebesgue measure on $\mathbb{R}^{M}$.
We show that the function $f$ on the right hand side is an integrable
majorant for all $\II_{A_{n}}(t)\,\mu_{n}\left(\left[-1,1\right]^{M}\backslash\left[-\tau,\tau\right]\right),n\in\mathbb{N}$:
\begin{align*}
\int_{\mathbb{R}^{M}}f(t)\textup{d}t & =\boldsymbol{\lambda}^{M}\left(A_{1}\right)\mu_{1}\left(\left[-1,1\right]^{M}\backslash\left[-\tau,\tau\right]\right)+\sum_{k=1}^{\infty}\boldsymbol{\lambda}^{M}\left(A_{k+1}\backslash A_{k}\right)\mu_{k}\left(\left[-1,1\right]^{M}\backslash\left[-\tau,\tau\right]\right).
\end{align*}
Each summand in the series above can be bounded above by
\begin{align*}
\boldsymbol{\lambda}^{M}\left(A_{k+1}\backslash A_{k}\right)\mu_{k}\left(\left[-1,1\right]^{M}\backslash\left[-\tau,\tau\right]\right) & \leq O\left(\left(\sqrt{k+1}-\sqrt{k}\right)\left(\sqrt{k+1}\right)^{M-1}\right)\mu_{k}\left(\left[-1,1\right]^{M}\backslash\left[-\tau,\tau\right]\right),
\end{align*}
which is summable in $k$ by (\ref{eq:summability}).

As $\left|\varphi_{\mathcal{N}\left(0,I_{M}\right)}\right|$ is integrable
as well, we have thus found an integrable majorant for $\II_{A_{n}}(t)\left|\varphi_{n}(t)-\varphi_{\mathcal{N}\left(0,I_{M}\right)}(t)\right|$.
By Theorem \ref{thm:fast_conv}, $\left|\varphi_{n}(t)-\varphi_{\mathcal{N}\left(0,I_{M}\right)}(t)\right|\xrightarrow[n\rightarrow\infty]{}0$
pointwise, so we conclude that the integral of\\
$\II_{A_{n}}(t)\left|\varphi_{n}(t)-\varphi_{\mathcal{N}\left(0,I_{M}\right)}(t)\right|$
over $\mathbb{R}^{M}$ converges to 0 as $n\rightarrow\infty$.

We proceed with the integrand over the set $B_{n}$:
\begin{align}
\II_{B_{n}}(t)\left|\varphi_{n}(t)\right| & \leq\II_{B_{n}}(t)\int_{\left[-1,1\right]^{M}}\prod_{\lambda}\left|\varphi_{\mathcal{R}\left(m_{\lambda}\right)}\left(\frac{t}{\sqrt{n_{\lambda}}}\right)\right|^{n_{\lambda}}\mu_{n}(\textup{d}m)\nonumber \\
 & \leq\II_{B_{n}}(t)\int_{\left[-1,1\right]^{M}}\left(\theta(m)\right)^{n}\mu_{n}(\textup{d}m),\label{eq:upper_bound_int_3}
\end{align}
where the existence of
\[
\theta(m)=\max_{t\in B_{n},\lambda=1,\ldots,M}\left|\varphi_{\mathcal{R}\left(m_{\lambda}\right)}\left(t_{\lambda}\right)\right|<1
\]
is a consequence of Lemma \ref{lem:char_fn}. We continue with the
calculation of an upper bound using the constant $\tau$ from (\ref{eq:summability}):
\[
(\ref{eq:upper_bound_int_3})\leq\II_{B_{n}}(t)\int_{\left[-\tau,\tau\right]}\left(\theta(m)\right)^{n}\,\mu_{n}(\textup{d}m)+\II_{B_{n}}(t)\,\mu_{n}\left(\left[-1,1\right]^{M}\backslash\left[-\tau,\tau\right]\right).
\]
On the interval $\left[-\tau,\tau\right]$, $\theta$ is bounded away
from 1:
\[
s:=\sup_{m\in\left[-\tau,\tau\right]}\theta(m)<1.
\]
This leads us to our final upper bound
\begin{align*}
\II_{B_{n}}(t)\left|\varphi_{n}(t)\right| & \leq\II_{B_{n}}(t)\left(s^{n}+\mu_{n}\left(\left[-1,1\right]^{M}\backslash\left[-\tau,\tau\right]\right)\right).
\end{align*}
For the last expression, we can construct an integrable majorant in
the same manner as for the second summand in (\ref{eq:upper_bound_int_2}).
This concludes the proof of Theorem \ref{thm:LLT}.

\section{\label{sec:Models}Examples of de Finetti Voting Models}

We discuss two classes of voting models featured prominently in the
literature that are de Finetti voting models according to Definition
\ref{def:de_Finetti_model}. As mentioned in the Introduction, models
of voting behaviour are often adapted from statistical physics. Two
categories are the collective bias model (CBM) and the Curie-Weiss
model (CWM). The single-group CBM first appeared in Straffin's work
on voting systems \cite{Straffin} with a uniform de Finetti measure.
It was later generalised in the works \cite{SRL,KL1}. The CWM in
its single-group version has a long history in statistical physics
and was first defined by Husimi \cite{Husimi} and Temperley \cite{Temperley}.
It has seen a multitude of applications to the social sciences (see
\cite{BD} for the first instance the CWM was used in such a context).
The multi-group version of the CWM was defined independently in \cite{CG}
and \cite{BRS}. Subsequently, the model has received a lot of attention
from other authors (see \cite{FC,Collet,FM,LS,KT1,KT2,KT3,KLSS}).

In the CWM, the voters tend to align with each other and there is
no external influence on them. In the CBM, the voters do not care
about others' opinions, but there is some central influence (such
as cultural or religious institutions) that induces correlation between
the voters. In this sense, the two models are opposite. However, mathematically,
they have more in common than a cursory glance would suggest. A multi-group
CBM is defined by the voting measure
\begin{equation}
\IP_{\textup{CBM}}\left(X_{11}=x_{11},\ldots,X_{Mn_{M}}=x_{Mn_{M}}\right):=\int_{\left[-1,1\right]^{M}}P_{m}^{\otimes n}\left(x_{11},\ldots,x_{Mn_{M}}\right)\mu\left(\textup{d}m\right)\label{eq:CBM}
\end{equation}
for all voting configurations $\left(x_{11},\ldots,x_{Mn_{M}}\right)$.
In the integral above, $\mu$ is a symmetric probability measure on
$\left[-1,1\right]^{M}$ and recall the probability measure $P_{m}^{\otimes n}$
from Notation \ref{not:prod_measure}. By (\ref{eq:CBM}), the CBM
is a de Finetti voting model according to Definition \ref{def:de_Finetti_model}
with a constant sequence of de Finetti measures $\mu_{n}:=\mu$, $n\in\IN$.

The variable $m$ in the integral (\ref{eq:CBM}) measures the prevalent
bias in the population. A positive bias means the voters are more
likely to be in favour of the proposal. For a fixed bias $m$, the
voters are conditionally independent in their decision. In a CBM,
\begin{equation}
\left(\frac{S_{n,1}}{n_{1}},\ldots,\frac{S_{n,M}}{n_{M}}\right)\xrightarrow[n\rightarrow\infty]{\textup{d}}\mu.\label{eq:CBM_lim}
\end{equation}
The above limit theorem provides useful information on the large population
behaviour of the CBM with the single exception of the (degenerate)
case $\mu=\delta_{0}$, in which all voters are independent and we
have the central limit theorem
\begin{equation}
\left(\frac{S_{n,1}}{\sqrt{n_{1}}},\ldots,\frac{S_{n,M}}{\sqrt{n_{M}}}\right)\xrightarrow[n\rightarrow\infty]{\textup{d}}\mathcal{N}\left(0,I_{M}\right).\label{eq:CBM_deg_lim}
\end{equation}
See \cite{KT_CBM} for these results.

While a CBM is a de Finetti voting model, it does not satisfy the
assumption of decreasing social cohesion embodied by the weak convergence
of $\mu_{n}$ to $\delta_{0}$, with the exception of $\mu=\delta_{0}$.
For more on single-group CBMs in the context of voting theory, see
\cite{SRL}, and for an in-depth discussion of multi-group versions
of the CBM, see \cite{KT_CBM}.

The single-group CWM is defined for all voting configurations $\left(x_{1},\ldots,x_{n}\right)$
by
\[
\IP_{\textup{CWM}}\left(X_{1}=x_{1},\ldots,X_{n}=x_{n}\right):=Z_{\beta,n}^{-1}\exp\left(\frac{\beta}{2n}\left(\sum_{i=1}^{n}x_{i}\right)^{2}\right),
\]
where $Z_{\beta,n}$ is a normalisation constant which depends on
$\beta$ and $n$. The parameter $\beta\geq0$ is the inverse temperature
in the physical context of the CWM as a model of ferromagnetism. As
a voting model, $\beta$ measures the degree of influence the voters
exert over each other. As we see, the most probable voting configurations
are those with unanimous votes in favour of or against the proposal.
However, there are only two of these extreme configurations, whereas
there is a multitude of low probability configurations with roughly
equal numbers of votes for and against. This is the `conflict between
energy and entropy'. Which one of these pseudo forces dominates depends
on the regime the model is in, determined by the value of the parameter
$\beta\geq0$. Using the Hubbard-Stratonovich transformation, the
CWM can also be expressed in the following way (see e.g.\!  Chapter
2 of \cite{FV2017} or Theorem 5.6 in the monograph \cite{MM}):
\begin{equation}
\IP_{\textup{CWM}}\left(X_{1}=x_{1},\ldots,X_{n}=x_{n}\right)=\left(Z_{\beta,n}'\right)^{-1}\int_{\left[-1,1\right]}P_{m}^{\otimes n}\left(x_{11},\ldots,x_{Mn_{M}}\right)\frac{\exp\left(-nF\left(m\right)\right)}{1-m^{2}}\,\textup{d}m,\label{eq:CWM_Finetti}
\end{equation}
for all voting configurations $\left(x_{1},\ldots,x_{n}\right)$.
\[
Z_{\beta,n}'=\int_{\left[-1,1\right]}\frac{\exp\left(-nF\left(m\right)\right)}{1-m^{2}}\,\textup{d}m
\]
is a normalising constant and the function $F$ above is given by
$F\left(m\right):=\frac{1}{2}\left(\frac{1}{\beta}\left(\frac{1}{2}\ln\frac{1+m}{1-m}\right)^{2}+\ln\left(1-m^{2}\right)\right),m\in\left[-1,1\right]$.
The representation (\ref{eq:CWM_Finetti}) of the CWM is called `de
Finetti representation'. We note that for fixed $n$ we can define
\begin{equation}
\mu_{n}:=\left(Z_{\beta,n}'\right)^{-1}\frac{\exp\left(-nF\left(m\right)\right)}{1-m^{2}}\,\II_{\left[-1,1\right]}\,\boldsymbol{\lambda},\label{eq:CWM_Finetti_measure}
\end{equation}
where $\boldsymbol{\lambda}$ is the Lebesgue measure, and the CWM
is a CBM. However, this is only true for fixed $n$. Whereas in (\ref{eq:CBM})
the de Finetti measure $\mu$ is static, i.e.\! independent of $n$,
in (\ref{eq:CWM_Finetti}) we see that there is a different de Finetti
measure $\mu_{n}$ for each model of size $n$. Hence, taken as a
class of models with any $n\in\IN$, the CWM is not a special case
of the CBM.

The CWM has three distinct `regimes' in which the voting margin
$S_{n}$ behaves differently. For $\beta<1$, the so called `high
temperature regime', $S_{n}/\sqrt{n}$ converges in distribution
to a centred normal distribution. The sequence of de Finetti measures
defined in (\ref{eq:CWM_Finetti_measure}) converges weakly to the
Dirac measure at 0. Moreover, it has a concentration property: for
any $\delta>0$, there are constants $C,D>0$ such that
\begin{equation}
\mu_{n}\left(\left[-1,1\right]\backslash\left[-\delta,\delta\right]\right)<C\exp\left(-Dn\right),\quad n\in\IN.\label{eq:exp_concentr_CW}
\end{equation}

This exponential decay is key to analysing the asymptotic behaviour
of the CWM. In the `critical regime', $\beta=1$, a limit theorem
holds for $S_{n}/n^{3/4}$ with a limiting distribution which is not
normal. Instead, the distribution has a density proportional to $\exp\left(-x^{4}/12\right),x\in\IR$.
In the `low temperature regime', $\beta>1$, $S_{n}/n$ converges
to the convex combination of two Dirac measures, $1/2\left(\delta_{-m(\beta)}+\delta_{m(\beta)}\right)$,
where $0<m(\beta)<1$ is the positive solution of the Curie-Weiss
equation $x=\tanh(\beta x)$.

The multi-group version of the CWM defined in (\ref{eq:CWM_Finetti})
is given by a positive semi-definite coupling matrix $J=\left(J_{\lambda\nu}\right)_{\lambda,\nu=1,\ldots,M}\in\IR^{M\times M}$
and the Gibbs measure
\begin{equation}
\IP_{\textup{CWM}}\left(X_{11}=x_{11},\ldots,X_{Mn_{M}}=x_{Mn_{M}}\right):=Z_{J,n}^{-1}\exp\left(\frac{1}{2}\sum_{\lambda,\nu=1}^{M}\frac{J_{\lambda\nu}}{\sqrt{n_{\lambda}n_{\nu}}}\sum_{i=1}^{n_{\lambda}}\sum_{j=1}^{n_{\nu}}x_{\lambda i}x_{\nu j}\right)\label{eq:multi-group_CW}
\end{equation}
for all voting configurations. The normalisation constant $Z_{J,n}$
depends on the group sizes $n_{\lambda}=n_{\lambda}(n)$ and the coupling
matrix $J$. The entries $J_{\lambda\nu}$ of the coupling matrix
give the strength of interaction between any pair of voters, one of
which belongs to group $\lambda$ and the other to group $\nu$. For
this probability measure, there is a de Finetti representation given
by
\begin{align*}
\IP_{\textup{CWM}}\left(X_{11}=x_{11},\ldots,X_{Mn_{M}}=x_{Mn_{M}}\right) & =\int_{\IR^{M}}P_{\bar{m}}^{\otimes n}\left(x_{11},\ldots,x_{Mn_{M}}\right)\exp\left(-nF_{J,n}\left(x\right)\right)\textup{{d}}x.
\end{align*}
Above, the mapping $m\mapsto\bar{m}$ is $\bar{m}:=\left(\tanh m_{\lambda}\right)_{\lambda=1,\ldots,M}$.
The CWM has a sequence of de Finetti measures given by
\begin{equation}
\mu_{n}:=\exp\left(-nF_{J,n}\right)\boldsymbol{\lambda}^{M},\quad n\in\IN,\label{eq:CWM_mu_n}
\end{equation}
with $\boldsymbol{\lambda}^{M}$ being the Lebesgue measure on $\IR^{M}$.
See \cite[display (12)]{KT3} for the form of $F_{J,n}$ and \cite[Section 2]{KT3}
for the three regimes of the model. An alternative representation
with compactly supported de Finetti measures can be obtained by a
change of variables in the integral in (\ref{eq:multi-group_CW})
setting $t:=\left(\tanh x_{\lambda}\right)_{\lambda=1,\ldots,M}$,
yielding de Finetti measures $\mu'_{n}$ whose support belongs to
$\left[-1,1\right]^{M}$. The multi-group CWM thus fits into the framework
outlined in Definition \ref{def:de_Finetti_model}.

The parameter space of the multi-group CWM is
\[
\Phi:=\left\{ \left.J\in\IR^{M\times M}\,\right|\,J\textup{ is positive definite or }J=\left(\beta\right)_{\lambda,\nu=1,\ldots,M}\textup{ for some }\beta\geq0\right\} .
\]
Analogously to the single-group model, the multi-group CWM has three
regimes. If $J$ is positive definite, the high temperature regime
is the region in $\Phi$ where $C:=\left(I_{M}-J\right)^{-1}$ is
positive definite. If $J=\left(\beta\right)_{\lambda,\nu=1,\ldots,M}$,
$\beta<1$ characterises the high temperature regime. In both cases,
this regime features weak interactions between voters. The sequence
$\left(\mu_{n}\right)_{n}$ given in (\ref{eq:CWM_mu_n}) exhibits
a critical convergence speed according to Definition \ref{def:conv_speed},
and the limiting distribution of $\left(S_{n,1}/\sqrt{n_{1}},\ldots,S_{n,M}/\sqrt{n_{M}}\right)$
is $\mathcal{N}\left(0,C\right)$ with a (generally) non-diagonal
covariance matrix $C$. As such, the limiting distribution of $\left(S_{n,1}/\sqrt{n_{1}},\ldots,S_{n,M}/\sqrt{n_{M}}\right)$
can be expressed as the convolution of two normal distributions:
\[
\left(\frac{S_{n,1}}{\sqrt{n_{1}}},\ldots,\frac{S_{n,M}}{\sqrt{n_{M}}}\right)\xrightarrow[n\rightarrow\infty]{\textup{d}}\mathcal{N}\left(0,I_{M}\right)\ast\mathcal{N}\left(0,\Sigma\right),
\]
where $\Sigma=\left(J^{-1}-I_{M}\right)^{-1}$ if the coupling matrix
$J$ is positive definite, and $\Sigma=\left(\beta/\left(1-\beta\right)\right)_{\lambda,\nu=1,\ldots,M}$
is singular if $J$ is a homogeneous matrix with all its entries equal
to some $0\leq\beta<1$. The critical regime of the CWM fits into
the subcritical regime according to Definition \ref{def:conv_speed}.
Lastly, the low temperature regime of the CWM features a sequence
$\left(\mu_{n}\right)_{n}$ defined in (\ref{eq:CWM_mu_n}) which
does \emph{not} converge weakly to $\delta_{0}$. See \cite[Sections 3 and 5.2]{KT3}
for more details.

In the CWM, all groups must be in the same regime of convergence speeds
given in Definition \ref{def:conv_speed}, unless there are clusters
of groups which are already independent for finite $n$ due to the
coupling constants being equal to 0 between groups belonging to different
clusters. This setup is realised by choosing the coupling matrix $J$
to be a block diagonal matrix. For CBMs, since the de Finetti sequence
$\left(\mu_{n}\right)_{n}$ is constant and equal to $\mu$ given
in (\ref{eq:CBM}), $\mu_{n}$ does not converge to $\delta_{0}$
at all, except in the degenerate case $\mu=\delta_{0}$ noted above.

No version of either the CBM or the CWM exhibits behaviour in line
with the empirical evidence provided in \cite{GKT,GKB} and discussed
in the Introduction. Depending on the de Finetti measure $\mu$, the
CBM produces expected per capita absolute voting margins for each
group $\lambda$, $\IE\left(\left|S_{n,\lambda}\right|/n_{\lambda}\right)$,
of order $1/\sqrt{n_{\lambda}}$ or $1$, as follows from the limit
theorems given in (\ref{eq:CBM_lim}) and (\ref{eq:CBM_deg_lim}),
with nothing in between. The CWM features three distinct regimes with
$\IE\left(\left|S_{n,\lambda}\right|/n_{\lambda}\right)$ of order
$1/\sqrt{n_{\lambda}}$, $1/n_{\lambda}^{1/4}$, and $1$, respectively.
(See \cite[Section 3]{KT3} for these results.)

We sum up the convergence behaviour of the CBM and the CWM in Table
\ref{tab:conv_speed}.

\begin{table}
\begin{centering}
\renewcommand{\arraystretch}{1.8}%
\begin{tabular}{|c|c|c|}
\hline 
CBM & CWM & De Finetti\tabularnewline
\hline 
\hline 
$\mu=\delta_{0}$ & $J=0$ \tablefootnote{Since $J$ is assumed to be either a positive definite or a homogeneous
matrix, this is indeed the only $J$ for which the CWM is in fast
convergence. If we allowed slightly more general $J$, such as block
diagonal matrices with one or more blocks equal to 0, then the groups
$\lambda$ for which $J_{\lambda\nu}=0$ for all $\nu=1,\ldots,M$
are in fast convergence, whereas other groups might fall into other
regimes. } & Fast convergence\tabularnewline
\hline 
-- & High temperature regime with $J\neq0$ & Critical convergence\tabularnewline
\hline 
-- & Critical regime & Slow convergence\tabularnewline
\hline 
$\mu\neq\delta_{0}$ & Low temperature regime & $\mu_{n}\nrightarrow\delta_{0}$\tabularnewline
\hline 
\end{tabular}
\par\end{centering}
\caption{\label{tab:conv_speed}De Finetti Regimes of the CBM and the CWM}

\end{table}

We next take a look at the asymptotic behaviour of the correlations
between two votes $\IE\left(X_{\lambda1}X_{\nu2}\right)$. By Proposition
\ref{prop:corr}, de Finetti voting models with an asymptotic loss
of social cohesion, represented by the weak convergence of the de
Finetti measures $\left(\mu_{n}\right)_{n}$ to $\delta_{0}$, feature
decaying correlations, i.e. $\IE\left(X_{\lambda1}X_{\nu2}\right)\xrightarrow[n\rightarrow\infty]{}0$.
Given the de Finetti measure $\mu$ of the CBM defined in (\ref{eq:CBM}),
we have the equality
\[
\IE\left(X_{\lambda1}X_{\nu2}\right)=\int_{\left[-1,1\right]^{M}}m_{\lambda}m_{\nu}\,\mu\left(\textup{d}m\right),\quad\lambda,\nu=1,\ldots,M.
\]
So correlations do not depend on the size of the population, and,
except for some special cases such as $\mu=\delta_{0}$, correlations
are non-zero in a CBM.

In the CWM, the asymptotic behaviour of the correlations $\IE\left(X_{\lambda1}X_{\nu2}\right)$
depends on the regime the model is in. We have (cf. \cite[Section 5.2]{KT3})
\begin{align*}
\IE\left(X_{\lambda1}X_{\nu2}\right) & =\begin{cases}
\Theta\left(\frac{1}{n}\right) & \textup{in the high temperature regime,}\\
\Theta\left(\frac{1}{\sqrt{n}}\right) & \textup{in the critical regime,}\\
\Theta\left(1\right) & \textup{in the low temperature regime.}
\end{cases}
\end{align*}
So we see that in a CWM the loss of social cohesion in the high temperature
and the critical regime is reflected in the decay of correlations
between pairs of votes. In the low temperature regime, where social
cohesion does not decay, neither do we have a decay of $\IE\left(X_{\lambda1}X_{\nu2}\right)$.

Finally, we provide a concrete example of a de Finetti voting model
in line with the cited empirical evidence:
\begin{example}
\label{exa:de_Finetti}Let $0.1\leq\alpha_{1},\ldots,\alpha_{M}\leq0.2$.
Let $\mu$ be a probability measure on the space $\left\{ -1,1\right\} ^{M}$
that is symmetric with respect to the origin. We set
\[
\mu_{n}:=\sum_{x\in\left\{ -1,1\right\} ^{M}}\mu\left\{ x\right\} \,\delta_{\left(\frac{1}{n_{1}^{\alpha_{1}}},\ldots,\frac{1}{n_{M}^{\alpha_{M}}}\right)\circ x},\quad n\in\IN,
\]
and let $\left(\IP_{n}\right)_{n}$ be the de Finetti voting model
with de Finetti sequence $\left(\mu_{n}\right)_{n}$. This model falls
into the contraction pattern defined in Theorem \ref{thm:phase_transition}.
Since $0.1\leq\alpha_{1},\ldots,\alpha_{M}\leq0.2$, its convergence
speed is subcritical according to Definition \ref{def:conv_speed}.
By Theorem \ref{thm:phase_transition}, we have the limit theorem
\[
\left(\frac{S_{n,1}}{n_{1}^{1-\alpha_{1}}},\ldots,\frac{S_{n,M}}{n_{M}^{1-\alpha_{M}}}\right)\xrightarrow[n\rightarrow\infty]{\textup{d}}\mu.
\]
Corollary \ref{cor:voting_margins} states that the expected per capita
voting margins exhibit the following behaviour:
\[
\IE\left(\frac{\left|S_{n,\lambda}\right|}{n_{\lambda}}\right)=\Theta\left(n_{\lambda}^{-\alpha_{\lambda}}\right),\quad\lambda=1,\ldots,M.
\]
Thus, choosing the parameters $\alpha_{\lambda}$, $\lambda=1,\ldots,M$,
allows us to obtain a good fit of the model to the data in \cite{GKT,GKB}.
Note that any possible correlation structure between the groups can
be realised by choosing the probability measure $\mu$ appropriately.
Hence, this choice can make the model fit the empirical correlations
between the voting margins $S_{\lambda}$ and $S_{\nu}$ of different
groups $\lambda$ and $\nu$. Due to both the discrete support and
the symmetry of $\mu$ with respect to the origin, this is a fairly
simple statistical model of binary voting. If desired, the introduction
of further symmetries in $\mu$ can lead to a reduction of the number
of parameters. See \cite{KT_CW_opt_weights} for a paper about optimal
voting weights in a two-tier voting system under a CWM with a reduced
set of parameters.
\end{example}

\section{\label{sec:Conclusion}Conclusion}

We discussed existing multi-group voting models such as the CBM and
the CWM. None of these models exhibits voting behaviour in line with
empirical evidence presented in \cite{GKT,GKB}. Therefore, we introduced
a general framework of probabilistic voting models defined by a de
Finetti representation. When the de Finetti measures converge weakly
to the Dirac measure at the origin, these models satisfy limit theorems
for the vector of group voting margins. The convergence speed of the
de Finetti measures falls into one of three regimes:
\begin{itemize}
\item If the convergence speed is fast (supercritical), which means faster
than $1/\sqrt{n}$, where $n$ is the population, then we obtain for
the normalised group voting margins a universal limiting distribution
which is multivariate normal with independent standard normal entries.
\item If the convergence speed is critical, which means of order $1/\sqrt{n}$,
and the de Finetti measures are a contraction of some underlying probability
measure $\mu$, then we obtain as the limiting distribution of the
normalised group voting margins a rescaled version of $\mu$ plus
an additive multivariate normal noise with independent entries.
\item Slow (subcritical) convergence speeds require normalisation of the
group voting margins by sequences which go to infinity faster than
$\sqrt{n}$, the normalisation factor in the other two regimes. The
limiting distribution of the normalised group voting margins is the
underlying probability measure $\mu$.
\end{itemize}
The subcritical regime can be used to model real-world voting behaviour
by assuming a convergence speed of order $n^{-\alpha}$, with $0.1\leq\alpha\leq0.2$,
such as in Example \ref{exa:de_Finetti}. By adjusting $\alpha$ to
the data, we can specify voting models in terms of some underlying
probability distribution $\mu$. This distribution can then be estimated
to obtain a statistical model of the voting data.

\section*{Acknowledgements}

This research was supported by a Secihti (formerly Conahcyt) postdoctoral
fellowship and through a SNII fellowship (candidate level) of the
author.

\end{document}